\documentclass[english]{elsarticle}
\usepackage[T1]{fontenc}
\usepackage[latin9]{inputenc}
\usepackage{geometry}
\usepackage{verbatim}
\usepackage{amsmath}
\usepackage{float}
\usepackage{amsbsy}
\usepackage{amssymb}
\usepackage{stackrel}
\usepackage{graphicx}
\usepackage{esint}
\usepackage{bm}
\usepackage{algpseudocode}
\usepackage{algorithm}
\usepackage{courier}

\usepackage{amsfonts} 

\usepackage{textcomp}
\usepackage{tabularx}
\usepackage{lineno}

\usepackage[table,xcdraw]{xcolor}
\usepackage{adjustbox,lipsum}
\usepackage{subcaption}

\usepackage{xcolor}

\newcommand{\zero}{\bm{0}}

\newcommand{\M}{\mathbf{M}}
\newcommand{\K}{\mathbf{K}}
\newcommand{\C}{\mathbf{C}}

\newcommand{\Ktg}{\mathbf{K^t}}
\newcommand{\Ktgr}{\mathbf{\tilde{K}^t}}

\newcommand{\Ktgrnl}{\mathbf{\tilde{K}^{(nl)}}}

\newcommand{\ktgr}{\tilde{\text{K}}^\text{t}}

\newcommand{\q}{\mathbf{q}}
\newcommand{\qd}{\mathbf{\dot{q}}}
\newcommand{\qdd}{\mathbf{\Ddot{q}}}

\newcommand{\qs}{\text{q}}

\newcommand{\fI}{\mathbf{f}}
\newcommand{\fIe}{\mathbf{f}_{e}}
\newcommand{\fIej}[1]{\textbf{f}_{#1}}

\newcommand{\KIej}[1]{\textbf{K}_{#1}}

\newcommand{\fe}{\mathbf{f}^{(\text{e})}}
\newcommand{\Ktge}{\mathbf{K}^{\mathbf{t}}_{e}}

\newcommand{\fIs}{\text{f}}

\newcommand{\Ko}{\mathbf{K}^{(1)}}
\newcommand{\Ktw}{\mathbf{K}^{(2)}}
\newcommand{\Kth}{\mathbf{K}^{(3)}}

\newcommand{\ko}{\text{K}^{(1)}}
\newcommand{\ktw}{\text{K}^{(2)}}
\newcommand{\kth}{\text{K}^{(3)}}

\newcommand{\Kor}{\mathbf{\tilde{K}}^{(1)}}
\newcommand{\Ktwr}{\mathbf{\tilde{K}}^{(2)}}
\newcommand{\Kthr}{\mathbf{\tilde{K}}^{(3)}}

\newcommand{\Korr}{\mathbf{\tilde{\tilde{K}}}^{(1)}}
\newcommand{\Ktwrr}{\mathbf{\tilde{\tilde{K}}}^{(2)}}
\newcommand{\Kthrr}{\mathbf{\tilde{\tilde{K}}}^{(3)}}

\newcommand{\fIr}{\mathbf{\tilde{f}}}
\newcommand{\fEr}{\mathbf{\tilde{f}^{\text{(e)}}}}

\newcommand{\kor}{\tilde{\text{K}}^{(1)}}
\newcommand{\ktwr}{\tilde{\text{K}}^{(2)}}
\newcommand{\kthr}{\tilde{\text{K}}^{(3)}}

\newcommand{\ktwri}{\overline{\text{K}}^{(2)}}
\newcommand{\kthri}{\overline{\text{K}}^{(3)}}

\newcommand{\Mr}{\mathbf{\tilde{M}}}
\newcommand{\Cr}{\mathbf{\tilde{C}}}

\newcommand{\V}{\mathbf{V}}
\newcommand{\VT}{\mathbf{V^T}}

\newcommand{\W}{\mathbf{W}}

\newcommand{\Ve}{\mathbf{V}_{e}}
\newcommand{\VeT}{\mathbf{V}^{\textbf{T}}_{e}}

\newcommand{\VejT}[1]{{\textbf{V}^{\mathbf{T}}_{#1}}}

\newcommand{\qj}[1]{\q^{{(#1)}}}
\newcommand{\qje}[2]{\q^{(#1)}_{#2}}

\newcommand{\et}{e}

\newcommand{\R}{\mathbb{R}}
\newcommand{\qr}{\bm{\eta}}
\newcommand{\qrd}{\mathbf{\bm{\dot{\eta}}}}
\newcommand{\qrdd}{\mathbf{\bm{\ddot{\eta}}}}

\newcommand{\qrs}{\eta}

\newcommand{\qrdds}{\ddot{\eta}}

\newcommand{\eigv}{\bm{\phi}}

\newcommand{\vv}[1]{\mathbf{v^{(#1)}}}

\newcommand{\Pder}[2]{\frac{\partial#1}{\partial #2}}

\newcommand{\Pderm}[3]{\frac{\partial^2 #1}{\partial #3 \partial #2}}

\newcommand{\tensCon}[2]{\cdot_{{#1}{#2}}}

\newcommand{\Els}{\mathbf{E}}

\newcommand{\Elsr}{\mathbf{\tilde{E}}}

\newcommand{\nel}{N_e}

\newcommand{\nelr}{\tilde{N}_e}

\newcommand{\ndofel}{n_e}

\newcommand{\el}{e}

\newcommand{\xiel}{\xi_e}

\newcommand{\xielvec}{\bm{\xi}}

\newcommand{\Nt}{N_s}

\newcommand{\Gecsw}{\mathbf{G}}
\newcommand{\becsw}{\mathbf{b}}

\newcommand{\gammaj}[1]{\bm{\gamma}^{(#1)}}

\newcommand{\delb}{\bm{\delta}}
\newcommand{\Rot}{\mathbf{U}}
\newcommand{\RotT}{\mathbf{U^T}}
\newcommand{\qrr}{\bm{\zeta}}
\newcommand{\qrrd}{\mathbf{\bm{\dot{\zeta}}}}
\newcommand{\qrrdd}{\mathbf{\bm{\ddot{\zeta}}}}

\newcommand{\sscr}[1]{\textsuperscript{\textcolor{blue}{#1}}}

\newcommand{\vm}{\bm{\phi}}
\newcommand{\vmT}{\bm{\phi}^T}

\def\onedot{$\mathsurround0pt\ldotp$}
\def\cddot{% two dots stacked vertically
  \mathbin{\vcenter{\baselineskip.67ex
    \hbox{\onedot}\hbox{\onedot}}%
  }}%
\def\cdddot#1{% three dots 
  \mathbin{\vcenter{\baselineskip.67ex
    \hbox{\onedot}\hbox{\onedot}\hbox{\onedot}%
  }}%
}

\newcommand{\eps}{\epsilon}



\usepackage{babel}

\journal{Computer Methods in Applied Mechanics and Engineering}

\begin{document}

\begin{frontmatter}{}

%\title{Efficient Reduced Order Model Construction through Hyperreduced Mesh for Random Vibration Analysis of Structural Components}

%\title{Accelerating Galerkin Reduced Order Model Tensors Identification through Hyperreduction}

\title{Accelerating Construction of Non-Intrusive Nonlinear Structural Dynamics Reduced Order Models through Hyperreduction}

\author{Alexander Saccani\corref{cor1}}
\ead{asaccani@ethz.ch}

\author{Paolo Tiso}

\address{Institute for Mechanical Systems,\\ ETH Z\"urich, \\ Leonhardstrasse 21, 8092 Z\"urich, Switzerland}

\begin{abstract}
We present a novel technique to significantly reduce the offline cost associated to non-intrusive nonlinear tensors identification in reduced order models (ROMs) of geometrically nonlinear, finite elements (FE)-discretized structural dynamics problems. The ROM is obtained by Galerkin-projection of the governing equations on a reduction basis (RB) of Vibration Modes (VMs) and Static Modal Derivatives (SMDs), resulting in reduced internal forces that are cubic polynomial in the reduced coordinates. The unknown coefficients of the nonlinear tensors associated with this polynomial representation are identified using a modified version of Enhanced Enforced Displacement (EED) method which leverages Energy Conserving Sampling and Weighting (ECSW) as hyperreduction technique for efficiency improvement.  Specifically, ECSW is employed to accelerate the evaluations of the nonlinear reduced tangent stiffness matrix that are required within EED. Simulation-free training sets of forces for ECSW are obtained from displacements corresponding to quasi-random samples of a nonlinear second order static displacement manifold. The proposed approach is beneficial for the investigation of the dynamic response of structures subjected to acoustic loading, where multiple VMs must be added in the RB, resulting in expensive nonlinear tensor identification. Superiority of the novel method over standard EED is demonstrated on FE models of a shallow curved clamped panel and of a nine-bay aeronautical reinforced panel modelled, using the commercial finite element program Abaqus.
\end{abstract}


\end{frontmatter}{}

%accelerating the eed method for ROM construction through hyperreduction 

%\pagebreak{}
%\linenumbers

%%%%%%%%%%%%%%%%%%%%%%%%%%%%%%%%%%%%%%%%%%%%%%%%%%%%%%%%%
%% Include additional sections %%%%%%%%%%%%%%%%%%%%%%%%%%
%%%%%%%%%%%%%%%%%%%%%%%%%%%%%%%%%%%%%%%%%%%%%%%%%%%%%%%%%
%\linenumbers
%% Introduction 
\section{Introduction} \label{sec:Introduction}
The investigation of the response of aerospace panels subjected to acoustic loading has been of great interest for the aerospace community in the past decades.
Experimental \cite{Hollkamp2018} and computational studies \cite{Gordon2011,Spottswood2008,Przekop2007,przekop2006nonlinear,mignolet2013review}  have been performed with the end goal of predicting fatigue life  and guarantee in service structural integrity.
With the advent of the \textit{Finite Element} (FE) method and the steady increase in computational power of modern computers, numerical investigations have become increasingly appealing compared to traditional experiments, which typically require costly testing facilities and instrumentation.\\
However, vibration assessment of real-life components using FEM is generally computationally demanding, if not intractable.  
This is due to the complex geometries that require large meshes, the necessity for simulations over extended time spans to obtain statistically relevant data, and the nonlinearities inherent in the \textit{Equations of Motion} (EOMs).
In fact, the pursue of light-design structures leads to large in service structural displacements whose behavior can be accurately modelled only by assuming finite rotations in the FE formulation.
As a result, \textit{geometric nonlinearities} are introduced in the FE model within the \textit{Total Lagrangian} formulation, or with the adoption of \textit{von K\'arm\'an} strains in thin walled structures undergoing mild deformations.  \par
\textit{Reduced Order Models} (ROMs) have been successfully employed to speed up the dynamic analysis of thin walled structures. 
The core principle behind model reduction is the replacement of solutions of the starting FE model (referred to as \textit{High Fidelity Model} (HFM)) with ROM counterparts which require very affordable computational effort.
%GP
Among different approaches, \textit{Galerkin Projection} (GP) is arguably the most classical and established technique for model reduction of structural dynamics problems, enjoying extensive application in academia as well as in an industrial setting. 
In GP ROMs, dimensionality reduction is achieved by restricting the HFM solution space to a small linear subspace spanned by few vectors constituting a \textit{reduction basis} (RB). \par
The choice of the RB dictates accuracy and efficiency of GP ROMs.
Different approaches, generally classifiable in  \textit{data-driven} and \textit{pyhsics-based} strategies, can be followed for RB construction.
%BASIS SELECTION
%data driven
Data-driven strategies rely on the availability of HFM solutions from which relevant vector components are usually extracted using \textit{Proper Orthogonal Decomposition} (POD) \cite{krysl2001dimensional,kerschen2007physical,lieu2005pod} and stacked together to form the RB. 
This comes at the large cost of running HFM simulations. 
%Physics based approach
On the other hand, physics-based approaches are simulation free and harness information from the physics modelled in the equations to construct the RB.
Among these techniques, \textit{Modal Analysis}  \cite{geradin2015mechanical} represents the cornerstone for model reduction of linear structural systems and acts as the foundation for extending reduction to nonlinear systems.
As such, the RB for geometrically nonlinear structures can be constructed by complementing the RB of the underlying linearized system with additional vectors that capture the most prominent nonlinear effects. 
To this end, \textit{Modal Derivatives} vectors (MDs) \cite{Idelsohn1985,idelsohn1985load} have been effectively utilized for model reduction of a large array of structures, both for transient and forced response \cite{Marconi2020,Marconi2021,Morteza,Sombroek2018,Wu2016,Rutzmoser2017,tiso2011reduction}.
 \par 
%Need for hyperreduction
Galerkin Projection usually enables a significant reduction in the number of degrees of freedom of the model, from hundreds of thousands, if not millions, to usually far less than a hundred. 
However, for nonlinear systems, the speed up in time integration of the ROM with respect to the HFM is generally small \cite{Farhat2014,Rutzmoser2017,Jain2018,An2008,rutzmoserThesis}, since another computational bottleneck arises:
the construction of the reduced internal force vector and of its Jacobian.
%In fact, these operations still retain the high dimensionality of the HFM, since the reduced forces and the reduced tangent stiffness matrix are computed by projecting the internal forces and tangent stiffness matrix of the HFM onto the RB. \par
For FE structures with geometric nonlinearities and linear elastic material, this bottleneck can be easily addressed by exploiting the polynomial nature of the internal forces \cite{Mignolet2008, mignolet2013review}. 
In an exact formulation, the reduced forces and their Jacobian are expressed as third order polynomials in the reduced coordinates and are be cheaply evaluated during time integration through tensor contraction operations. 
Conversely, the computation of tensors associated to reduced forces is performed once and for all offline, i.e. during model construction prior to time integration.  \par
To this purpose, different techniques, broadly classifiable in \textit{direct methods} and \textit{indirect methods}, can be employed \cite{mignolet2013review}.
The direct methods, were the first to be used \cite{nashThesis, shi1996finite, almroth1978automatic} and rely on the knowledge of the decomposition of the internal forces into their linear, quadratic and cubic components.
This separation is generally not provided by commercial FE packages like ABAQUS, NASTRAN or ANSYS, and as such, had limited the use of tensorial formulation to ROM constructed from \textit{ad-hoc} written FE codes. 
This limitation was overcome in \cite{Muravyov2003,McEwan2001,Hollkamp2008} with the introduction of \textit{indirect} methods, where tensors are identified using the FE code as a black-box, returning forces to displacements inputs and displacements to forces inputs.
The most popular indirect methods for tensor identification are the \textit{Implicit Condensation} (IC)\cite{McEwan2001}, \textit{Implicit Condensation and Expansion} (ICE)\cite{Hollkamp2008} and \textit{Enforced Displacement} (ED)\cite{Muravyov2003}. \par
A major limitation of indirect schemes is their high offline computational cost for large ROMs featuring many vectors in the RB \cite{Perez2014}.
Models for acoustic load predictions of aerospace components fall usually into this category, since the high modal density and the large excitation bandwidth of the load force the analyst to include a large number of vectors in the RB \cite{Perez2014}.
The \textit{Enhanced Enforced Displacement} (EED) was proposed in this context \cite{Perez2014}. The EED cuts the computational cost of ED tensor identification through the use of the tangent stiffness matrix instead of the internal force vector only.
 \\
%why we still need to accelerate tensor construction
Even with EED, tensor construction remains the computational bottleneck in the construction of large \textit{physics-based} ROMs, potentially constituting a major issue when multiple realizations are needed.
This scenario can be envisioned in geometry optimization problems or in uncertainty quantification in presence of geometrical defects \cite{Wang2018,Perez2011}. \par
With the end goal of reducing overall model construction time, one could  consider abandoning the tensorial formulation and resorting to \textit{hyperreduction} techniques \cite{farhat2020}, where the reduced order forces and their Jacobians are replaced with a quickly computable approximation, during time integration. 
Within these methods we mention \textit{Discrete Empirical Interpolation} (DEIM) \cite{chaturantabut2010nonlinear} and \textit{Energy Conserving Sampling and Weighting} (ECSW) \cite{Farhat2014,Farhat2015}. 
However, both DEIM and ECSW reduced forces approximation is constructed from training forces that are generally obtained from HFM simulations and could compromise ROM offline construction time.
 Simulation-free strategies have been proposed in  \cite{Rutzmoser2017, Jain2018} for ECSW only.
Specifically, ECSW was trained in \cite{Rutzmoser2017,rutzmoserThesis} with nonlinear static solutions, whereas in \cite{Jain2018} with displacements obtained from a linear run lifted on a second order nonlinear manifold. \par
Physics-based ROMs, equipped with ECSW hyperreduction, constructed with simulation-free training, allow potentially for a significant reduction in model construction time when compared to ROMs in the tensorial formulation, as reported in  \cite{rutzmoserThesis}. 
The drawback of ECSW-ROMs lies, however, in their intrusive nature: during solution of the ROMs equations, internal element forces and jacobians are computed with the FE code.
This fact does not constitute a problem if the numerical integration scheme and ROM construction are embedded in the FE package, while it becomes a major limitation when the FE formulation is not accessible.
This is usually the case when ROMs equipped with ECSW hyperreduction rely on commercial FE packages: the communication cost between ROM solver and FE code severely undermines online ROM performance as recently reported in \cite{Trainotti2024}.
%what we propose
\par
Based on these premises, with the ultimate goal of fully non-intrusive ROMs, we propose a new method to accelerate ROM tensor identification from commercial FE packages. Our approach is based on EED identification procedure and ECSW trained with the second order manifold as in \cite{Jain2018}. The idea is to cheaply construct a hyperreduced mesh to use for a fast evaluation of the tangent stiffness in the EED scheme.
Our method shares similarities with the recent work proposed in \cite{Kim2023}.
In that work, an ECSW reduced mesh was constructed following the same approach proposed in \cite{Jain2018} and used to speed up tensor construction of a parametric ROM within a direct identification scheme.
In this work we propose a similar strategy based on ECSW for efficient tensor identification, coupled however with an indirect identification scheme.
This great distinctive feature of the method here proposed from the one presented in \cite{Kim2023}, allows for efficient tensor construction using commercial FE packages, rendering the method applicable in an industrial setting.  \par
%content
The content of this paper is as follows:
In section 2 we recall the basics of Galerkin ROMs focusing on RBs composed of MDs, tensor idenfication based on EED and ECSW hyperreduction strategy; In section 3 we illustrate the core of the proposed approach for tensor identification; 
Eventually, in section 4, we provide the application of the methodology to two case study: a rectangular curved panel and an aerospace multi-bay stiffened panel subjected to acoustic loading.

% preliminaries
\section{Preliminaries} 

\label{sec:ROM_MDs}

\subsection{Governing Equations for FE Structures}
\label{sec:HFMeqs}
The equations of motion for FE models of thin-walled structures with geometric nonlinearities and linear elastic material read
\begin{equation}
    \M \qdd + \C\qd + \fI (\q) = \fe(t),    
    \label{eq:eomF}
\end{equation}
where $\q \in \R^{n}$ is the vector of nodal displacements, $\M, \C \in \R^{n \times n}$ are the mass and viscous damping matrices, $\fI(\q), \fe \in \R^n$ are, respectively, the vectors of internal and external forces. The tangent stiffness matrix $\Ktg \in \R^{n \times n}$ is given as
\begin{equation}
    \Ktg = \Pder{\fI}{\q}.
\end{equation} 
In these models, the kinematic strain-displacements relations are based on the \textit{Green-Lagrange} strain tensor or on the approximate \textit{von K\'arm\'an} strains \cite{crisfield1991}. 
This leads to an internal force vector $\fI(\q)$ that is an exact cubic polynomial function of the nodal displacements $\q$ \cite{Mignolet2008}, which can be be written in tensorial and Einstein's notation as
\footnote{Unless differently specified, we adopt in the following Einstein's summation convention over repeated indices.}
\begin{equation}
    \fI(\q) = \Ko \cdot \q + \Ktw \cddot (\q \otimes \q) + \Kth \cdddot\   (\q \otimes \q \otimes \q )
    \iff
    \fIs_{i} = \ko_{ij}\qs_{j} + \ktw_{ijk}\qs_{j}\qs_{k} + \kth_{ijkl} \qs_{j}\qs_{k}\qs_{l},
    \label{eq:tensDecF}
\end{equation}
where $\Ko \in \mathbb{R}^{n\times n}$ is the linear stiffness matrix, $\Ktw \in \mathbb{R}^{n\times n \times n}$ and $\Kth \in \R^{n \times n \times n \times n}$ are respectively the tensors associated to the quadratic and cubic forces, "$\otimes$" is the dyadic product operator, and '$\ \cdot\ $','$\ \cddot\ $','$\ \cdddot\ \ $' denote respectively, single, double and triple contraction. 
 The tensors $\Ktw$ and $\Kth$ are commonly referred to as \textit{nonlinear stiffness operators}, while $\Ko$ is the linear stiffness.
The tensorial formulation of internal forces in Eq.\eqref{eq:tensDecF} can be derived by partial differentiation of a quartic order elastic potential function.
As a result, the FE tensors are symmetric with respect to all axis, i.e. they are invariant to permutations of their indices.
For example, $\ktw_{ijk} = \ktw_{ikj} = \ktw_{jik} = \ktw_{jki} = \ktw_{kij} = \ktw_{kji}$.
As will be shown in the sequel, this decomposition of internal forces is essential for efficient model reduction.
In the following, Eq. \eqref{eq:eomF} is referred to as HFM. Note that the method we propose could also be applied when other formulations, as for instance co-rotational, are adopted. In this case, the form \eqref{eq:tensDecF} is an approximation of the actual non-polynomial reduced order forces.

%Galerkin Projection
\subsection{Galerkin Projection}
In GP ROMs, model reduction is achieved by restricting the HFM solution to a small linear subspace of the HFM solution space.
In particular, it is enforced that
\begin{equation}
    \q(t) \approx \V \qr(t),
    \label{eq:redBas}
\end{equation}
where $\V \in \R^{n\times m}$ (with $m \ll n $) is the RB, and $\qr \in \R^{m}$ is the vector of reduced coordinates.
After substitution of Eq. \eqref{eq:redBas} into Eq. \eqref{eq:eomF}, the derived equations are projected onto the reduction basis $\V$, obtaining
\begin{equation}
    \VT\M\V\qrdd + \VT\C\V\qrd + \VT\fI(\V\qr) = \VT\fe,
\label{eq:redEqsGP}
\end{equation}
where we dropped in the notation the dependence on the reduced coordinates on time.
The ROM equations can be equivalently derived by applying the principle of \textit{Virtual Work} to the HFM with the  constraint in  Eq.\eqref{eq:redBas} by requiring that the residual $\mathbf{r}$ obtained inserting \eqref{eq:redBas} into \eqref{eq:eomF} is orthogonal to $\V$, i.e. $\VT \mathbf{r} = \zero$. 
As a result, the Galerkin ROM is said to preserve the \textit{Lagrangian Structure} \cite{Farhat2015} associated to the HFM. \\
Choice of the RB impacts on accuracy and efficiency of the derived ROM: 
a good RB should span as accurately as possible the HFM solution with a small number of vectors.
%basis selection
\subsection{Basis construction with Static Modal Derivatives}
\label{sec:MDs}
\subsubsection{Reduction Basis for Linear Systems: Vibration Modes}
Physics-based approaches represent an excellent choice for RB selection, since they are computationally cheap, as opposed to data-driven methods.
For structural dynamics, the RB is usually constructed starting from a set of \textit{Vibration Modes} (VMs) of the underlying linearized system.
These vectors are computed by solving the following eigenvalue problem:
\begin{equation}
    (-\M \omega_i^2 + \Ko) \eigv_i = \zero,\ \ \ \text{for}\ \ 
    i = 1,...,n,
    \label{eq:eigLin}
\end{equation}
where $\eigv_i$ is the $i$\textsuperscript{th} VM and $\omega_i$ its corresponding angular frequency.
Strategy for VMs selection depends on the application. 
For unforced transient analysis, the relevant VMs are the ones with lowest frequencies, since high frequency modes are generally heavily damped.
For forced problems, however, VMs should  be selected both based on the \textit{Modal Participation Factor} (MPF) of the load \cite{geradin2015mechanical}, and on its frequency content: all the VMs with non zero MPF, whose eigenfrequency is within (or close to) the excitation frequency bandwidth of the load, should be included in the RB. 
\subsubsection{Reduction Basis for Nonlinear Systems: Static Modal Derivatives}
Using RBs of only VMs for model reduction of nonlinear structures leads generally to poor results. 
In fact, some additional vectors, modelling nonlinear effects, should be included in the RB.
In thin-walled structures, geometric nonlinearities are due to coupling between out-and in-plane displacements: as structural deformation occurs along out plane modes (which are generally the ones excited by the load) axial stretching arises \cite{tiso20214}.\\ 
Modal Derivatives vectors, first presented in \cite{Idelsohn1985,idelsohn1985load} and successively employed for reduction of a wide array of structures  \cite{Marconi2020,Marconi2021,Morteza,Sombroek2018,Wu2016,Andersson2023}, proved to well capture the bending-stretching coupling.
These vectors were defined in \cite{Idelsohn1985} as the sensitivities of the eigenvalue problem in Eq. \eqref{eq:eigLin} to a change of the linearization point in the directions of the VMs. 
However, only their simplified version, known as \textit{Static Modal Derivatives} (SMDs)\cite{Jain2017}, have been mostly used for practical model reduction problems. 
Therefore, we recall in this section the definition and properties of SMDs only, referring the reader to \cite{Idelsohn1985,Weeger2016} for a detailed discussion on MDs.\\
The SMD $\bm{\theta}_{ij} \in \R^n$, related to VMs $\bm{\phi}_i,\bm{\phi}_j$, is defined as
\begin{equation}
    \Ko\bm{\theta}_{ij} = - \left. \Pderm{\fI(\q)}{\eps_i}{\eps_j} \right|_{\eps_i,\eps_j = 0} ,
    \label{eq:SMDdef}
\end{equation}
with $\q = \bm{\phi}_i\eps_i+\bm{\phi}_j\eps_j$.
When the decomposition of the internal forces in Eq. \eqref{eq:tensDecF} holds, the above definition for SMDs is equivalent to
\begin{equation}
    \Ko\bm{\theta}_{ij} = - \Ktw \cddot (\bm{\phi_i} \otimes \bm{\phi_j}),
    \label{eq:SMDdef1}
\end{equation}
and the SMDs share a symmetry property for which $\bm{\theta}_{ij} = \bm{\theta}_{ji}$ \cite{Jain2017}, since the FE tensors are symmetric.  
This allows to construct RB for the nonlinear system by complementing the RB for its linearized counterpart with only $({n_{\phi}}^2+{n_{\phi}})/2$ vectors,  instead of  ${n_{\phi}}^2$ vectors, being ${n_{\phi}}$ the number of retained VMs.
Notice that the number of SMDs grows almost quadratically with the number of VMs. 
This could lead to large RBs, resulting in inefficient ROMs. 
To mitigate this problem, a selection strategy for SMDs, based on a linear run, has been proposed in \cite{tiso2011optimal}. 
\par
The SMDs are computed by solving the linear problem in Eq.\eqref{eq:SMDdef}. The derivative term at r.h.s. can be retrieved non-intrusively exploiting the chain differentiation rule and by approximating the directional derivative of the tangent stiffness matrix with finite differences \cite{Morteza}.
For example, when a central difference scheme is adopted, the approximation reads
\begin{equation}
    \left. \Pderm{\fI(\q)}{\eps_i}{\eps_j} \right|_{\eps_i,\eps_j = 0} 
    =  \left.\Pder{\Ktg(\q)}{\epsilon_i}\right|_{\eps_i,\eps_j = 0} 
    \bm{\phi}_j
    \approx \frac{\Ktg{(h\bm{\phi}_i)}-\Ktg{(-h\bm{\phi}_i)}}{2h}\bm{\phi}_j,
    \label{eq:SMDsRHS}
\end{equation}
where $h$ is a user defined scalar perturbation parameter.
\subsubsection{Static Quadratic Manifold}
\label{sec:SQM}
In addition to the classic projection approach, SMDs have been used in a second order manifold for reduction of thin walled structures \cite{Jain2017}. Loosely speaking, the idea behind the quadratic manifold is to constrain SMDs amplitudes, to the amplitudes of VMs.
More specifically, starting with $n_{\phi}$ VMs, the HFM solution is enforced to lay on the quadratic constraint that reads
\begin{equation}   \q = \bm{\Gamma}(\bm{\gamma}) = \sum_{i=1}^{n_{\phi}} \gamma_i \bm{\phi}_i + \frac{1}{2}\sum_{i=1}^{n_{\phi}}\sum_{j=1}^{n_{\phi}} \gamma_i \gamma_j\bm{\theta}_{ij},
\label{eq:QM}
\end{equation}
where $\gamma_i$ is the amplitude of VM $\bm{\phi}_i$. In the following, we refer to this manifold as the \textit{Static Quadratic Manifold} (SQM).\par
For a von-K\'arm\'an beam, the SQM is equivalent to the \textit{Static Condensation} \cite{Jain2017}, a technique in which the axial degrees of freedom are statically constrained to the transverse ones. 
Even if model reduction with the SQM proved to work effectively for some flat structures \cite{Jain2017,Rutzmoser2017QM}, its application to slightly curved structures can lead to misleading predictions \cite{Vizzaccaro2021}.
In fact, in more recent work in \cite{Vizzaccaro2021} and \cite{Veraszto2020} it was shown that the SQM is generally not invariant, unless the spectral gap between the frequency of the VMs retained in the expansion and the ones of the other VMs is large.
For this reason, we  use in the following the more robust approach where the SMDs are complementing a RB of VMs.
%TENSORIAL FORMULATION
\subsection{Tensorial Formulation}
\label{sec:tensForm}
Despite model reduction, integrating the ROM equations as in \eqref{eq:redEqsGP} typically incurs a computational cost comparable to that of the HFM.
The reason behind this limited speed-up is the expensive construction of the reduced force vector and of its Jacobian during Newton-Raphson iterations within the solution scheme.
In fact, both these two quantities are obtained by projection of their counterpart associated to the HFM. 
As such, internal forces and tangent stiffness matrices must be first constructed using the FE model with \textit{element-level} and \textit{assembly} operations, whose computational cost scales with the size of the HFM. \par
This large cost can be avoided by
exploiting the polynomial formulation of the internal forces in \eqref{eq:tensDecF}.
In particular, upon substitution of Eq \eqref{eq:redBas} into Eq. \eqref{eq:tensDecF}, the term $\VT \fI(\V \qr)$ in Eq.\eqref{eq:redEqsGP} can be written as
\begin{equation}
    \fIr(\qr) \triangleq \VT \fI(\V \qr) =  \Kor \cdot \qr + \Ktwr \cddot (\qr \otimes \qr) + \Kthr \cdddot\   (\qr \otimes \qr \otimes \qr )
    \iff
   \tilde{\fIs}_{i} = \kor_{ij}\qrs_{j} + \ktwr_{ijk}\qrs_{j}\qrs_{k} + \kthr_{ijkl} \qrs_{j}\qrs_{k}\qrs_{l},
    \label{eq:redForcTens}
\end{equation}
where the reduced tensors $\Kor \in \R^{m \times m}$, $\Ktwr \in \R^{m \times m \times m}$ and $\Kthr \in \R^{m \times m \times m \times m}$ are defined as
\begin{subequations}
\begin{align}
    \Kor &= \VT \Ko \V \hspace{5cm} 
     &\kor_{ij} = V_{li}\ko_{lk} V_{kj}\\ 
    \Ktwr &= \left( \VT \Ktw \cdot \V \right)\tensCon{2}{1} \V \hspace{3cm}
    \iff 
    &\ktwr_{ijk} = V_{li} \ktw_{lrp} V_{pk} V_{rj}\\ 
    \Kthr &= \left( \left( \VT \Kth \cdot \V \right) \tensCon{3}{1} \V \right) \tensCon{2}{1} \V
    &\kthr_{ijkl} = V_{pi} \kth_{pqrs} V_{qj} V_{rk} V_{sl}. 
\end{align}
\label{eq:tensProj}
\end{subequations}
In the above formulas and in the following, "$\tensCon{i}{j}$" denotes tensor contraction of the $i$\textsuperscript{th} dimension of the first tensor with the $j$\textsuperscript{th} dimension of the second tensor (e.g. for $\mathbf{A} \in \R^{k\times k \times k}$ and $\mathbf{B} \in \R^{k\times k}$, $\mathbf{A}\tensCon{2}{1} \mathbf{B}$ is equivalent to $\text{A}_{ijk}\text{B}_{jl}$), and we have adopted Einstein's summation convention over repeated indices. \\
In similar fashion,  the polynomial expression for the reduced tangent stiffness matrix, defined as \begin{equation}
   \Ktgr \triangleq  \frac{\partial }{\partial \qr}\VT \fI(\V \qr) = \VT \Ktg (\V \qr) \V
   \label{eq:redTgstiffDef}
\end{equation}
writes
%\begin{equation}
%    \tilde{\text{K}}^t_{ij} = \kor_{ij} + \text{K}_{ijp}^{(t2)} \eta_{p} + \text{K}_{ijpq}^{(t3)}\eta_p\eta_q 
 %   \label{eq:redTgSt}
%\end{equation}
\begin{equation}
    \Ktgr = \Kor + \mathbf{\tilde{K}^{(t2)}} \qr + \mathbf{\tilde{K}^{(t3)}}\cddot \left( \qr \otimes \qr \right)  \hspace{1cm} \iff \hspace{1cm}
    \ktgr_{ij} = \kor_{ij} + \tilde{\text{K}}_{ijp}^{(t2)} \qrs_{p} + \tilde{\text{K}}_{ijp}^{(t3)} \qrs_{p} \qrs_{q},
    \label{eq:redTgSt}
\end{equation}
where
\begin{subequations}
\begin{align}
    \tilde{\text{K}}_{ijp}^{(t2)} &= \ktwr_{ijp}+\ktwr_{ipj} \\
 \tilde{\text{K}}_{ijpq}^{(t3)} &= \kthr_{ijpq}+\kthr_{ipqj}+\kthr_{ipjq}.
\end{align} 
\label{eq:redTgSt1}
\end{subequations}
With the expression for the reduced forces $\fIr(\qr)$ provided in  Eq.\eqref{eq:redForcTens}, the ROM equations in \eqref{eq:redEqsGP} compactly write
\begin{equation}
    \Mr \qrdd + \Cr \qrd + \fIr(\qr) = \fEr(t), \\
    \label{eq:ROMtens}
\end{equation}
where 
\begin{equation}
 %   \fIr = \Kor \cdot \qr + \Ktwr \cddot\   (\qr \otimes \qr) + \Kthr \cdddot\   (\qr \otimes \qr \otimes \qr ), \quad
    \fEr(t) \triangleq \VT \fe(t), 
    \quad \Mr \triangleq \VT \M \V,
    \quad \Cr \triangleq \VT \C \V.
\end{equation}
In this formulation, the reduced order forces and their Jacobian can be efficiently computed during time integration by contracting the tensors $\Kor$, $\Ktwr$, $\Kthr$, $\mathbf{\tilde{K}^{(t2)}}$ and $\mathbf{\tilde{K}^{(t3)}}$ onto the vector of reduced coordinates $\qr$, according to Eqs. \eqref{eq:redForcTens} and 
\eqref{eq:redTgSt}. 
%In this formulation, if the expression for the tensors $\Kor$, $\Ktwr$, $\Kthr$, $\mathbf{\tilde{K}^{(t2)}}$ and $\mathbf{\tilde{K}^{(t3)}}$$\Kor$, $\Ktwr$, $\Kthr$, $\mathbf{\tilde{K}^{(t2)}}$ and $\mathbf{\tilde{K}^{(t3)}}$ are known, the evaluation of the reduced forces and of their Jacobian can be efficiently performed during time integration with tensor-vector multiplications as in Eqs.\eqref{eq:redForcTens} and \eqref{eq:redTgSt}.\\
\\
Reduced order tensors inherit the same symmetry property of their full counterpart, as described in section \ref{sec:HFMeqs}. 
Exploiting symmetries enables efficient storage, faster online reduced forces evaluations and provides a metric to assess accuracy of the identified tensors \cite{rutzmoserThesis}.\par
\subsubsection{Methods for Tensors Identification}
In order to exploit the substantial advantages of the tensorial formulation for time integration, tensors must be correctly identified in the model construction process.
Strategies for tensors identification can be classified in two different categories: \textit{direct methods} and \textit{indirect methods} \cite{mignolet2013review}.
Direct methods were the first 
 to be proposed in \cite{nashThesis, shi1996finite, almroth1978automatic}.
 Within these intrusive techniques, reduced tensors are computed by direct projection of the full tensors on the RB, performing the operations in Eq. \eqref{eq:tensProj}.
 As such, direct methods require the knowledge of the full order tensors $\Ko$, $\Ktw$ and $\Kth$, limiting  ROM applicability to FE codes that offer this output. \par
Conversely, indirect methods are non-intrusive, in the sense that tensors are constructed from commonly provided FE outputs, i.e. from nodal forces to imposed displacements or nodal displacement to imposed forces. 
This distinctive feature of indirect methods makes them appealing for the solution of industrial problems, that are usually tackled with commercial FE programs. 
Indirect methods can be further divided in Implicit Condenstation (IC) methods \cite{McEwan2001,Hollkamp2008} and ED methods \cite{Muravyov2003,Perez2014}. 
The former (IC) are used to construct ROMs where the membrane inertia is neglected and membrane contributions to generalized forces are implicitly modelled in the identified tensors. 
In fact, the RB of IC-ROMs is usually composed of only bending VMs.
On the other hand, ED methods are employed in ROMs where the nonlinear stretching effects are explicitly modelled with additional vectors in the RB, such as \textit{Dual Modes} \cite{Mignolet2008}, MDs and SMDs. 
%TENSOR IDENTIFICATION
%direct and indirect methods
\subsection{Enhanced Enforced Displacements Method}
\label{sec:EED}
The ED method, also known as \textit{STiffness Evaluation Procedure} method (STEP), was first put forward in \cite{Muravyov2003}. 
This scheme for tensor identification can be broken down in two subsequent steps. 
In the first step, 
the FE solver is used to compute nonlinear nodal forces to a set of prescribed displacements in the direction of the vectors of the RB.
In the second step, linear algebraic systems of equations are set up with the previously computed forces and solved for the unknown tensor entries.
For $m$ vectors in the RB, $(m^3-3m^2+2m)/6$ nonlinear static problems must be solved to fully identify the tensors. 
For this reason, tensor identification can become prohibitively expensive, as ROM size increases.\par
The EED was proposed in \cite{Perez2014} to mitigate the high computational cost of ED.
The Enhanced Enforced Displacements leverages the output of the tangent stiffness matrix to  reduce the total number of nonlinear static solutions required for tensor identification.
Using EED, tensors can be fully identified from only $(m^2+5m)/2$ nonlinear static solutions. In this section, we present a short review of the EED method, as our new approach  for tensor identification here presented is heavily based on it.
\subsubsection{Alternative Expression for Tensorial Forces}
The idea behind EED is the identification of tensors entries exploiting  Eqs.\eqref{eq:redTgSt}\eqref{eq:redTgSt1} that relate the reduced tangent stiffness matrix to reduced coordinates and tensor coefficients.
Before delving into the identification algorithm, we present an alternative formulation of Eqs.\eqref{eq:redForcTens},\eqref{eq:redTgSt} and \eqref{eq:redForcTens1}, on which tensor identification is based. 
The internal forces are re-written as
\begin{equation}
     \tilde{f}_i = (\text{V}^T\fIs)_{i} = \sum_{j=1}^{m}\kor_{ij}\qrs_{j} + \sum_{j=1}^m\sum_{k= j}^m\ktwri_{ijk}\qrs_{j}\qrs_{k} + \sum_{j=1}^m \sum_{k = j}^m\sum_{l=k}^m\kthri_{ijkl} \qrs_{j}\qrs_{k}\qrs_{l}, \hspace{1cm} i\in \{1,...,m\},
    \label{eq:redForcTens1}
\end{equation}
while the reduced tangent stiffness as
\begin{equation}
    \begin{split}
    \ktgr_{ij} = \Pder{ }{\qrs_j}(\text{V}^T\fIs)_{i} = \kor_{ij} + &\sum_{p=j}^m \ktwri_{ijp}\qrs_p + \sum_{p=1}^{j} \ktwri_{ipj}\qrs_p 
    + \sum_{p=j}^m\sum_{q=p}^m \kthri_{ijpq} \qrs_p \qrs_q + \\
    + & \sum_{p=1}^m\sum_{q= p}^{j} \kthri_{ipqj}\qrs_{p}\qrs_q+\sum_{p=i}^{j} \sum_{q=j}^m \kthri_{ipjq}\qrs_p \qrs_q, \hspace{1cm} i,j\in \{1,...,m\},
    \end{split}
    \label{eq:tgStiffId}
\end{equation}
where we restrict the summations for $\ktwri_{ijk}$ to $j\leq k$ and for $\kthri_{ijkl}$ to $j\leq k \leq l$, and impose all other entries to be zero. 
In fact, many coefficients in Eq.\eqref{eq:redForcTens} multiply the same monomial of the reduced coordinates and thus, their coefficients could be summed together in a single coefficient (e.g. $\ktwr_{ij}\qrs_i\qrs_j + \ktwr_{ij}\qrs_i\qrs_j$ can be rewritten as $(\ktwr_{ij}+\ktwr_{ji})\qrs_i\qrs_j$ if $i \neq j$. 
This alternative reformulation of the reduced forces is essential to avoid indeterminacy in the identification process, since every monomial of the reduced coordinates is multiplied only by one unique coefficient. 
From here on, to avoid a cumbersome notation, we denote $\kthri$ with $\kthr$, and $\ktwri$ with $\ktwr$. 
The reader should, however, keep in mind that the identified tensors are sparse, containing a non zero unique entry for each  monomial in the reduced coordinates.
\subsubsection{Tensor Identification}
The first step in EED tensor identification consists of imposing to the FE model, nodal displacements $\q$ in the direction of the vectors of the RB. \ 
If we denote with $\textbf{v}^{(r)}$, the $r$\textsuperscript{th} column vector in the RB, for a displacement of the form $\q = \qrs_r \textbf{v}^{(r)}$, Eq.\eqref{eq:tgStiffId} can be simplified to three different cases (no Einstein's summation):
\begin{subequations}
    \begin{align}
    \ktgr_{ij} - \kor_{ij} &= \ktwr_{irj}\qrs_r + \kthr_{irrj}\qrs_r^2 \hspace{1cm} &\text{if} \quad r < j,\\
    \ktgr_{ij} - \kor_{ij} &= 2\ktwr_{irr}\qrs_r + 3\kthr_{irrr}\qrs_r^2 \hspace{1cm} &\text{if} \quad r = j,\\
    \ktgr_{ij} - \kor_{ij} &= \ktwr_{ijr}\qrs_r + \kthr_{ijrr}\qrs_r^2 \hspace{1cm} &\text{if} \quad r > j,
\end{align}
 \label{eq:identifEqs}
\end{subequations}
with $i,j,r\in\{1,...,m\}$.
The only unknowns in these equations are the coefficients of the nonlinear tensors $\Ktwr$ and $\Kthr$. 
In fact, the reduced coordinates $\qrs_r$ and the linear stiffness $\Kor$ are known\footnote{$\Kor$ is obtained by projecting the linear tangent stiffness onto the RB}, while the reduced tangent stiffness $\Ktgr$ is computed with Eq.\eqref{eq:redTgstiffDef} from the full tangent stiffness  $\Ktg(\V\qr)$ returned by the FE solver. 
In all the three cases, Eqs.\eqref{eq:identifEqs} contain two unknowns each. 
As such, tensor entries are retrieved by imposing two different displacements in the direction of each vector $\textbf{v}^{(r)}$ in the RB ($2m$ in total), and by solving a set of two by two linear systems of equations. 
In this way, all the coefficients of the form $\ktw_{ijl}$, $\kthr_{ijjl}$, $\kthr_{iljj}$ and $\kthr_{ijjj}$ are identified.\par
Coefficients in the cubic force tensor, corresponding to triplets of different indices (i.e. $\kthr_{irsj}$ with $i,r,s,j \in\{1,..,m\},\  r<s<j$), are identified in a subsequent step.
The expression of the reduced tangent stiffness for displacements of the form $\q = \qrs_s \vv{s} + \qrs_r \vv{r}$ writes (no Einstein's summation):
\begin{equation}
    \ktgr_{ij} - \kor_{ij} = \ktwr_{irj}\qrs_r + \ktwr_{isj}\qrs_s + \kthr_{irsj}\qrs_r\qrs_s + \kthr_{irrj}\qrs_r^2 + \kthr_{issj}\qrs_r^2  
    \hspace{1cm} 
    \text{with} \hspace{1cm} r<s<j,
     \label{eq:identifEqs1}
\end{equation}
As can be easily seen, the only unknown in this equation is the sought $\kthr_{irsj}$, since all the other tensor coefficients have been previously identified. 
Exploiting this, tensors identification is  completed  by imposing to the FE model additional $m(m-2)/2$ displacements of the form $\q = \qrs_s \vv{s} + \qrs_r \vv{r}$ for $s\neq r$ and by setting up and solving, for each of them, the linear Eq.\eqref{eq:identifEqs1}.
\subsection{Energy Conserving Sampling and Weighting}
\label{sec:ECSW}
\subsubsection{Approximate Evaluation of Reduced Forces}
Energy Conserving Sampling and Weighting hyperreduction \cite{Farhat2014,Farhat2015} is an alternative to tensorial formulation for efficient computation of reduced order forces.
An approximation of these forces and of their Jacobian is obtained by evaluating forces and tangent stiffness for a small subset of elements of the FE mesh,
referred to as \textit{reduced mesh}.\\
If we denote with $\Els$ the set of elements in the FE mesh, and $\Elsr \subset \Els$ the set of elements in the reduced mesh, the ECSW approximation writes
\begin{subequations}
\begin{align}
   \VT \fI(\V\qr) = \sum_{\et \in \Els}^{|\Els|} \VeT \fIe(\Ve \qr) \ \approx \   
   \sum_{\et \in \Elsr}^{|\Elsr|} \xiel \VeT \fIe(\Ve \qr),\\
     \VT \Ktg(\V\qr) \V = \sum_{\et \in \Els}^{|\Els|} \VeT \Ktge(\Ve \qr) \Ve \ \approx \   
   \sum_{\et \in \Elsr}^{|\Elsr|} \xiel \VeT \Ktge(\Ve \qr) \Ve,
   \label{eq:ECSWred}
\end{align}
\end{subequations}
where $\fIe \in \R^{\ndofel}$ is the vector of nodal forces of element $\el$, $\Ve \in \R^{\ndofel\times m}$ is the restriction of the RB to the rows of nodal displacements of element $\el$, $\Ktge \in \R^{\ndofel\times\ndofel}$ is the tangent stiffness matrix of element $\el$ and $\xiel \in \R$ is the weight associated to element $\el$. We define $\nel = |\Els|$ and $\nelr = |\Elsr|.$
Loosely speaking, the contributions of missing elements to the reduced forces (which have the units of energy) is accounted for by weighting the contributions from the elements retained in the reduced mesh, ensuring that total energy is preserved.
The speed-up in reduced forces construction achievable with ECSW is proportional to the number of elements left out of the reduced mesh. 
In other words, the smaller $\Elsr$, the greater the computational gain.
\subsubsection{Computation of Reduced Mesh and Weights}
The reduced mesh, i.e. the elements in $\Elsr$ and their associated weights $\xiel$, is found by solving a minimization problem based on a set of training nodal displacements of the FE model $\{\qj{1},\ldots,\qj{\Nt}\}$. 
If we denote with $\qje{i}{\et}$ the nodal displacements of element $\et$ for snapshot $i$, we can define the quantities $\Gecsw \in \R^{\Nt \cdot m \times \nel} $ and $\becsw \in \R^{\Nt \cdot m}$ as
\begin{equation}
    \Gecsw = 
    \begin{bmatrix}
         \VejT{1} \fIej{1}(\qje{1}{1}) & \ldots & \VejT{\nel}\fIej{\nel}(\qje{1}{\nel})\\
         \vdots & \ddots & \vdots\\
         \VejT{1} \fIej{1}(\qje{\Nt}{1})  & \ldots & \VejT{\nel}\fIej{\nel}(\qje{\Nt}{\nel})
    \end{bmatrix}, \quad
    \becsw = 
    \begin{bmatrix}
         \sum_{\et = 1}^{\nel}\VejT{\et} \fIej{\et}(\qje{1}{\et}) \\ 
         \vdots \\
         \sum_{\et = 1}^{\nel} \VejT{\et} \fIej{\et}(\qje{\Nt}{\et}) 
    \end{bmatrix},
    \label{eq:ECSWGb}
\end{equation}
and the constrained minimization in $\xielvec \in \R^{\nel}$ problem as
\begin{equation}
    \xielvec^* = \text{argmin} \| \xielvec \| _{0}  \hspace{1cm} \text{s.t.} \hspace{1cm} \| \Gecsw\xielvec-\becsw \|_2 \leq \tau \| \becsw \|_2, \ \xi_i \geq 0\quad \text{for}\quad i\in \{1,\ldots\nel\},
    \label{eq:minProbl}
\end{equation}
where  $\|\bullet \|_{i}$ denotes the $i$\textsuperscript{th} norm of vector $\bullet$, and $0<\tau<1$ is a user defined tolerance.  \\
In words, the optimal solution $\xielvec^*$ is the one that minimizes the number of non-zero entries in $\xielvec$ while reproducing the total work of the internal forces in vector $\becsw$ with a relative error less then $\tau$.
The zero norm in the objective function definition in Eq.\eqref{eq:minProbl} makes the problem \textit{NP-hard} and thus computationally intensive to solve. 
As such, its optimal solution is replaced with a sub-optimal solution searched with greedy algorithms, such as the \textit{sparse Non-Negative Least Square} algorithm (sNNLS). 
The reader is referred to \cite{Farhat2014} for more details on the solution scheme. \par
In contrast to hyperreduction strategies based on interpolation, such as DEIM, ECSW preserves the \textit{Lagrangian Structure} of the internal forces model.
In other words, if the internal forces are conservative, their hyperreduced approximated version is conservative as well \cite{Farhat2015}.
Moreover, the non-negative constraint on the weights $\xi$ preserves the positive definiteness of the reduced tangent stiffness matrix. 
The combination of these two factors ensures that the hyperreduced model has the same numerical stability properties of the Galerkin exact ROM \cite{Farhat2015}.
%and $\xielvec \in \R^{\nel}$
%mention also DEIM

%method
\section{Efficient Tensors Construction using EED-ECSW} 
\label{sec:ID}
As opposed to ROMs with hyperreduction, ROMs in tensorial formulation, once constructed, do not retain any bound with the FE code.
This enables efficient integration, that would not be achievable with ROMs relying on hyperreduction approximation based on commercial FE packages. 
In this latter case, during integration, the time spent for communication between the ROM solver and the FE code increases significantly the overall solution time \cite{Trainotti2024}.
\par
A limitation of this approach lies in the high tensor construction time when many vectors in the RB are used, as is usually the case for ROMs employed for fatigue life assessment of thin-walled structures \cite{Perez2014, Gordon2011, Spottswood2008}.
In this section we propose a novel method to speed up tensor construction for projection ROMs with RB of VMs and MDs.
The presented technique is fully non-intrusive, from the computation of the RB to tensor construction, making it applicable in an industrial setting.
The method consists of accelerating EED for tensor identification, using ECSW trained with displacements from the SQM presented in section \ref{sec:SQM}.
In the following, we will refer to this method as \textit{EED-ECSW}.\par
%Tensor identification is performed in three subsequent steps: computation of the ECSW model, improved identification with EED and tra
%We present them in the following in two dedicated sections.\\
\subsection{Reduction Basis: Effect of Orthogonalization on Tensor Identification}
\label{sec:need4Trans}
Before delving into the proposed approach, we briefly discuss 
RB construction and its implication in the subsequent tensors identification. 
The ROM herein presented is based on a RB of VMs and SMDs  that, with the same notation introduced in section \ref{sec:MDs}, writes
\begin{equation}
    \V = \left[\bm{\phi}_1,\ldots,\bm{\phi}_{n_{\phi}},\bm{\theta}_{11},\ldots,\bm{\theta}_{n_{\phi}n_{\phi}}\right].
    \label{eq:physicalRB}
\end{equation} 
Selection strategies for SMDs can be potentially employed to reduce RB size.
Good practice for numerically stable ROMs is the orthogonalization of the RB prior to projection, using numerical orthogonalization schemes such as the \textit{Gram-Schmidt} algorithm. 
We denote with $\W$ the orthogonalized counterpart of $\V$ in Eq.\eqref{eq:physicalRB}.
A proper choice of $\W$ has to ensure good conditioning of the ROM equations, thus enabling numerically stable ROM integration \cite{rutzmoserThesis}.\par
In our experience, the displacements shapes of column vectors in $\W$ are generally irregular and less smooth as compared to those of vectors in $\V$.
Imposing displacements to the FE model in the direction of these vectors, within EED,  can potentially lead to poorly identified tensors, as their shape do not represent, in general, meaningful displacement fields for tensor identification.
We attribute these issues to the fact that the FE formulation for internal forces is, in many cases, close to, but not exactly a cubic polynomial. 
This is the case, for instance, for FE codes implementing the \textit{Updated-Lagrangian} formulation or the \textit{co-rotational} formulation for shell elements.
These deviations from the cubic formulation make the identified tensors dependent on the choice of imposed displacements, as recently suggested in \cite{Lin2023}.
In this work, we circumvent this hurdle by identifying the tensors using the physical RB $\V$ and by subsequently transform them to their counterpart corresponding to the orthogonalized RB $\W$, using transformation relations introduced in the sequel. 
\subsection{Identification based on EED-ECSW}
The core idea we put forward in this manuscript is to speed up EED tensor identification scheme (see section \ref{sec:EED}) by approximating the reduced tangent stiffness matrix with ECSW.
\label{sec:EED_ECSW}
Specifically, the left hand side of Eqs. \eqref{eq:identifEqs} and \eqref{eq:identifEqs1} can be approximated as:
\begin{equation}
    \Ktgrnl \triangleq \Ktgr - \Kor = \V^T(\Ktg-\Ko)\V \approx \sum_{\el \in \Elsr}^{\nelr} \xiel \VeT \left[\Ktge(\q = \Ve \qr) - \KIej{\et}\right] \Ve,
    \label{eq:approxECSW}
\end{equation}
where we remind the reader that $\KIej{\et}$ and $\Ktge$ are respectively the reduced linear stiffness matrix and reduced tangent stiffeness matrix of element $\et$.\\
The computational advantage over standard EED lies in the fact that the tangent stiffness matrices are computed only for the elements in the reduced mesh.
In this way, remarkable time savings are achieved in the operations performed by the FE code and in reading the entries of the outputted tangent stiffness matrices, which are generally fewer. 
\par
\subsection{Computation of the ECSW Model}
\label{sec:ECSWmodel}
Tensor identification based on EED-ECSW promises to be very efficient, as compared to EED, when a small reduced mesh is employed.
However, an additional overhead cost has to be sustained to find the ECSW reduced mesh and the associated weights.
As such, the overall efficiency of the proposed method, can be strongly affected by this cost, that must be kept as contained as possible.
For this reason, we could not afford to train ECSW with forces snapshots coming from computationally expensive HFM simulations, as is usually done in the literature.
Alternative ways for computationally affordable training were put forward in \cite{Rutzmoser2017,Jain2018}.
In \cite{Rutzmoser2017} training snapshots were collected from static solutions to pseudo-inertial loads, whereas in \cite{Jain2018}, they were obtained from the evaluation of nodal forces for displacements coming from a linear run, lifted with the SQM presented in section \ref{sec:MDs}.
This second strategy is more convenient when using RB with SMDs, since the SQM can be assembled from these vectors at low additional computational cost.
Hence, we decided to adopt a modified version of this last approach for ECSW training.
\subsubsection{Static Quadratic Manifold based Training for ECSW}
\label{sec:SQM_4train}
Our training is based on the SQM presented in Eq.\eqref{eq:QM}, constructed from the VMs retained in the RB.
Firstly, a set of nodal displacements of the FE model $\mathbf{S_{\q}}=\{\qj{1},\ldots,\qj{{\Nt}}\}$ is assembled by evaluating the SQM for a set of  VMs amplitudes $\mathbf{S_{\gamma}} = \{\gammaj{1},\ldots,\gammaj{\Nt}\}$ as
\begin{equation}
    \qj{i} = \bm{\Gamma}(\gammaj{i}). 
\end{equation}
The set of VMs amplitudes $\mathbf{S_{\gamma}}$ is obtained using \textit{Latin-Hypercube Sampling} \cite{helton2003latin}. 
This quasi-random sampling scheme allows for optimal coverage of the parameter space.
Bounds for the sampled variables are defined for each entry $\gamma_i$ of vector $\bm{\gamma}$, based on the maximum out of plane displacement associated to the corresponding VM.
Specifically, we impose that 
\begin{equation}
    \gamma_i \in [-\delta_i,+\delta_i] \hspace{1cm} \text{with} \hspace{1cm} \delta_i = \frac{\alpha}{\text{max}(\bm{\phi}_i)},
    \label{eq:scaling}
\end{equation} 
where $\bm{\phi}_i$ is the VM associated to $\gamma_i$, $\alpha$ is a positive user provided constant and '$\text{max}(\star)$' returns the maximum entry of vector '$\star$'.
In this way, for each displacement in $\mathbf{S_{\q}}$, we ensure that every component of the linear part of the SQM $\gamma_i \bm{\phi}_i$ (see Eq. \eqref{eq:QM}) does not exceed displacements of magnitude $\alpha$.
Enforcing this constraint is necessary to guarantee that the displacements returned from the manifold are physically meaningful: 
a good value of $\alpha$ should force the structure to deform in the nonlinear regime, at a displacement level similar to the one experienced in the dynamic response. 
In practical cases, $\alpha$ can be chosen to be a fraction of the thickness of the structure at hand.
In fact, if the VMs in the RB are out of plane dominated, the maximum displacement in vector $\bm{\phi}_i$ is likely to be an out of plane dof, allowing for a direct comparison with the thickness.
Notice that when shell elements are used, the vector $\bm{\phi}_i$ in Eq.\eqref{eq:scaling} should be restricted to displacements dofs, excluding rotational dofs, for which a bound is more difficult to be defined. \\ \par 
With the displacements in $\mathbf{S_q}$ available,
a set of unassembled nodal element forces $\mathbf{S_f}$ is computed using the FE model:
\begin{equation}
    \mathbf{S_f} = 
     \begin{Bmatrix}
        \begin{bmatrix}
            \fIej{1}(\qje{1}{1}) \\ \vdots \\ 
            \fIej{\nel}(\qje{1}{\nel})
        \end{bmatrix}, &
        \ldots, &
         \begin{bmatrix}
            \fIej{1}(\qje{\Nt}{1}) \\ \vdots \\ 
            \fIej{\nel}(\qje{\Nt}{\nel})
        \end{bmatrix}
    \end{Bmatrix},
    \label{eq:forceSet}
\end{equation}
where we adopt the same notation used in section \ref{sec:ECSW}, in which $\fIej{\et}$ is the vector of nodal forces of element $\et$, $\qje{j}{\et}$ is the nodal displacement vector of snapshot $j$ restricted to the dofs of element $\et$ and $\nel$ is  the total number of elements in the FE model.\\
The approximation in \eqref{eq:approxECSW} features the difference between the reduced tangent stiffness matrix and the reduced linear stiffness matrix.
This difference can be viewed as the jacobian of the nonlinear part of the reduced forces:
\begin{equation}
   \Ktgrnl = \V^T(\Ktg-\Ko)\V = \Pder{}{\qr}\VT\left( \fI(\V\qr)-\Ko\V\qr\right)
\end{equation}
Hence, in order to be consistent with Eq. \eqref{eq:approxECSW}, we need to train ECSW to reproduce
\begin{equation}
    \VT\left( \fI(\V\qr)-\Ko\V\qr\right),
\end{equation}
instead of simply $\VT \fI(\V\qr)$, as is usually done.
As such, we
subtract linear forces contributions to samples in $\mathbf{S_f}$ before assembling  $\Gecsw$ and $\becsw$ matrices, thus obtaining the following set of training forces
\begin{equation}
    \mathbf{S_{f_{nl}}} = 
     \begin{Bmatrix}
        \begin{bmatrix}
            \fIej{1}(\qje{1}{1}) - \KIej{1}\qje{1}{1}\\ \vdots \\ 
            \fIej{\nel}(\qje{1}{\nel}) - \KIej{\nel}\qje{1}{\nel}
        \end{bmatrix}, &
        \ldots, &
         \begin{bmatrix}
            \fIej{1}(\qje{\Nt}{1}) - \KIej{1}\qje{\Nt}{1}\\ \vdots \\ 
            \fIej{\nel}(\qje{\Nt}{\nel}) - \KIej{\nel}\qje{\Nt}{\nel}
        \end{bmatrix}
    \end{Bmatrix}.
    \label{eq:forceSetNl}
\end{equation}
where $\KIej{\et}$ is the linear stiffness of element $\et$.\\
\subsubsection{Computation of ECSW Reduced Mesh}
\label{sec:redMeshComp}
The force set defined in Eq.\eqref{eq:forceSetNl} is used to determine the reduced mesh and weights in the ECSW scheme presented in section \ref{sec:ECSW}. 
This set is divided into a training set and a validation set, with $N_t$ and $N_v$ samples, respectively. 
Two different pairs of $\Gecsw$ and $\becsw$ matrices are assembled from the forces in these two sets, $\Gecsw_t$, $\becsw_t$ and $\Gecsw_v$, $\becsw_v$, respectively, using the physical RB $\V$ defined in Eq.\eqref{eq:physicalRB}.
For example $\Gecsw_t$ and $\becsw_t$ used for training write
\begin{gather}
    \Gecsw_t = 
    \begin{bmatrix}
         \VejT{1} (  \fIej{1}(\qje{1}{1}) - \KIej{1}\qje{1}{1} ) & \ldots & \VejT{\nel} ( \fIej{\nel}(\qje{1}{\nel}) - \KIej{\nel}\qje{1}{\nel} )\\
         \vdots & \ddots & \vdots\\
         \VejT{1} (\fIej{1}(\qje{N_t}{1}) - \KIej{1}\qje{N_t}{1} )  & \ldots & 
         \VejT{\nel} ( \fIej{\nel}(\qje{N_t}{\nel}) - \KIej{\nel}\qje{N_t}{\nel} )
    \end{bmatrix},\\
    \becsw_t = 
    \begin{bmatrix}
         \sum_{\et = 1}^{\nel}\VejT{\et}(   \fIej{\et}(\qje{1}{\et}) - \KIej{\et}\qje{1}{\et}  )  \\ 
         \vdots \\
         \sum_{\et = 1}^{\nel} \VejT{\et} (   \fIej{\et}(\qje{N_t}{\et}) - \KIej{\et}\qje{N_t}{\et}  ) 
    \end{bmatrix}.
    \label{eq:ECSWGb}
\end{gather}
A reduced mesh and ECSW weights  are found by solving the minimization problem in Eq.\eqref{eq:minProbl} with the $\Gecsw$ and $\becsw$ matrices associated to the training set, for specified relative tolerance $\tau$.  
The solution algorithm that we adopt in this work is the sNNLS presented in \cite{Jain2018} (Algorithm 1).\\
The tolerance parameter $\tau$ is critical for effectiveness of the proposed approach. 
If $\tau$ is too small, the reduced mesh is larger than necessary and the speed-up of the subsequent tensor identification is undermined.
Conversely, if $\tau$ is too large, tensor identification can be inaccurate, leading to poor dynamic predictions of the resulting ROM.
Tolerance values used in the literature \cite{Farhat2015, Rutzmoser2017}, in the range $0.0001\leq\tau\leq0.01$  proved to be valid for the proposed approach, as reported in the result section \ref{sec:applications}. 
\par
Extrapolation performance of the ECSW model for reduced forces predictions on the SQM can be assessed with the validation set. 
Given the optimal weights $\xielvec^*$, a scalar relative error $\epsilon_{\text{ECSW}}$ can be computed as 
\begin{equation}
    \epsilon_{\text{ECSW}} = \frac{\| \Gecsw_{v} \xielvec^* - \becsw_{v} \|_2}{\| \becsw_{v} \|_2}.
    \label{eq:ecswErr}
\end{equation}
If the validation error is too large, one could think of increasing the number of training samples or alternatively of solving again the minimization problem for lower tolerance $\tau$.
In both cases, the number of elements in the reduced mesh is expected to increase, while the  the relative error $\epsilon_{\text{ECSW}}$ to decrease.\\
A summary of simulation-free ECSW model generation presented in this section is provided in Algorithm \ref{alg:redMesh}.
%algorithm for reduced mesh
\begin{algorithm}[H]
\caption{ECSW Reduced Mesh Computation}
\textbf{Input:} \texttt{FeModel}; RB $\V$; VMs $\bm{\Phi}$\sscr{a}; SMDs $\bm{\Theta}$\sscr{b}; bounding parameter $\alpha$; number of training and validation samples $N_t$,$N_v$; tolerance for ECSW $\tau$\\
\textbf{Output}:  ECSW reduced element set $\Elsr$; element weights $\xielvec^*$
\label{alg:redMesh}
\begin{algorithmic}[1] 
    \Statex \textit{Generation of random VMs amplitudes samples} 
    \State $\bm{\delta} \gets \texttt{boundVMs}(\bm{\Phi},\alpha)$\sscr{c}
    \State $N_s \gets N_t + N_v$
    \State $ \{\gammaj{1},\ldots,\gammaj{\Nt}\}
    \gets \texttt{LHS}(\bm{\delta},\Nt)$\sscr{d}
    \Statex \textit{SQM lifting} 
    \State Build function \texttt{SQM}: $\bm{\Gamma}(\bm{\gamma}) \gets \texttt{SQM}(\bm{\gamma})$\sscr{e}
    \For{$i \in \{1,2,\ldots,N_s\}$}
    \State $\qj{i} \gets \texttt{SQM}(\gammaj{i})$
    \EndFor
    \Statex \textit{Computation of training forces}
    \State \{$\KIej{1},...,\KIej{\nel}\} \gets \texttt{FeModel.ElementStiffLin()}$\sscr{f}
    \State $\mathbf{S_{f_{nl}}} \gets \texttt{tuple}(N_s)$\sscr{g} 
    \For{$i \in \{1,2,\ldots,N_s\}$}
    \State $ \{ \fIej{1},...,\fIej{\nel} \} \gets  \texttt{FeModel.ElementForces}(\qj{i})$\sscr{h}
    \State $\mathbf{S_{f_{nl}}}\{i\} \gets \texttt{subtractLinearForce}(\{ \fIej{1},...,\fIej{\nel} \},\{\KIej{1},...,\KIej{\nel}\},\qj{i}) $\sscr{i}
    \EndFor
    \Statex \textit{Train ECSW model}
    \State $[\mathbf{G}_{\bm{t}},\mathbf{b}_{\bm{t}},\mathbf{G}_{\bm{v}},\mathbf{b}_{\bm{v}}] \gets \texttt{assembleGb($\V,\mathbf{S_{f_{nl}}},N_t,N_v$)}$\sscr{l}
    \State $[\xielvec^*,\Elsr] \gets \texttt{SNNLS}(\mathbf{G}_{\bm{t}},\mathbf{b}_{\bm{t}},\tau)$\sscr{m}
    \Statex \textit{Validate ECSW model}
    \State $\epsilon_{\text{ECSW}} \gets \texttt{relErrECSW}(\mathbf{G}_{\bm{v}},\mathbf{b}_{\bm{v}},\xielvec^*)$\sscr{n}
\end{algorithmic}
\end{algorithm}
\vspace{-10pt}
\footnotesize
\noindent \sscr{a}$\bm{\Phi}$ array containing all the VMs in RB $\V$.\\
\sscr{b}$\bm{\Theta}$ array containing all the SMDs in RB $\V$.\\
\sscr{c}\texttt{boundVMs} function returns the bounds of Latin Hypercube according to Eq.\eqref{eq:scaling}.\\
\sscr{d}\texttt{LHS} function returns quasi-random samples of VMs amplitude vectors $\bm{\gamma}$.\\
\sscr{e} Assemble the \texttt{SQM} function from VMs and SMDs, according to Eq.\eqref{eq:QM}.\\
\sscr{f} \texttt{FeModel.ElementStiffLin}$()$ function returns the linear stiffnesses of all the elements in \texttt{FeModel}.\\
\sscr{g}\texttt{tuple}$(i)$ initializes a tuple array with $i$ elements.\\
\sscr{h}\texttt{FeModel.ElementForces}$(\qj{i})$ returns nodal forces for all the elements in \texttt{FeModel}, corresponding to displacements $\qj{i}$.\\
\sscr{i}\texttt{subtractLinearForce} function subtracts the linear force to total force as in Eq.\eqref{eq:forceSetNl}.\\
\sscr{l} \texttt{assembleGb} function assembles $\Gecsw$ and $\becsw$ arrays defined in Eq.\eqref{eq:ECSWGb}, both for ECSW training and validation.\\
\sscr{m}\texttt{SNNLS} function solves optimization problem in Eq.\eqref{eq:minProbl} with the sparse Non-Negative Least Square algorithm.\\
\sscr{n} \texttt{relErrECSW} function returns the validation error of the ECSW model, as in Eq.\eqref{eq:ecswErr}.\\
\normalsize
\subsubsection{Comparison with Training from SQM Lifted-linear Run}
The training here presented differs from the one in  \cite{Jain2017,Kim2023}, in which the modal displacements in set $\bm{S_{\gamma}}$ are coming from a linear modal solution to the same load used in the nonlinear analysis.
In fact, this approach can potentially lead to non-physical forces. 
As proved in \cite{Weeger2016}, the SQM in Eq.\eqref{eq:QM} is the second order Taylor expansion of the solution to the nonlinear static problem
\begin{equation}
    \fI(\q) = \sum_{i=1}^{n_{\phi}} \gamma_i \Ko \bm{\phi_i},
\end{equation}
where $\bm{\phi_i}$ is the $i$\textsuperscript{th} VM.
Displacements from linear runs can be potentially larger (or smaller) as compared to those of nonlinear runs and, when inserted into Eq.\eqref{eq:QM}, can return non-physical displacements, if the input is far-off the range of validity of the expansion.
With our new approach we enforce the limits of the modal coordinates, possibly preventing this problem.
On top of that, the response of a linear system can be far less rich in terms of modal content, as compared to the response of its nonlinear counterpart, since nonlinear phenomena such as internal resonances cannot be captured with a linear model.
For this reason, the training proposed in \cite{Jain2018} may not be rich enough in general.
Conversely, the choice of sampling the modal coordinates space with LHS allows to optimally cover a wide range of dynamic responses, where multiple modes are activated at the same time.

\subsection{Tensor Identification and Transformation to Orthogonalized RB}
\label{sec:tensRot}
With the computed ECSW reduced mesh, tensor identification is performed using EED, where we employ the approximation in Eq. \eqref{eq:approxECSW} for the construction of the reduced tangent stiffness.
In this way we can identify the ROM tensors $\Kor$, $\Ktwr$ and $\Kthr$ for the physical RB $\V$.
 These tensors need to be transformed to their counterpart associated with the orthogonalized RB $\W$, prior to ROM integration.
 In the following, we denote the transformed tensors with $\Korr$, $\Ktwrr$ and $\Kthrr$. \\
 Formulas to change between two different RB can be derived by applying  the principle of Virtual Work to the ROM in Eq.\eqref{eq:ROMtens} obtaining
 \begin{equation}
     \delb\qr^{\bm{T}} \left[ \Mr \qrdd + \Cr\qrd + \fIr(\qr)-\fEr(t) \right] = \bm{0} \hspace{1cm} \forall \ \delb\qr \in \R^m,
 \end{equation}
 where $\delb\qr$ is an infinitesimal variation of the reduced coordinates.
 Given the change of variables in reduced coordinates space
 \begin{equation}
     \qr = \Rot \qrr, \hspace{1cm} \Rot \in \R^{m \times m},
     \label{eq:rotation}
 \end{equation}
 the above condition re-writes
 \begin{equation}
     \delta \qrr^{\bm{T}} \left[ \Rot^{\bm{T}} \Mr \Rot \qrrdd + \Rot^{\bm{T}} \Cr \Rot \qrrd + \Rot^{\bm{T}} \fIr(\V\Rot\qrr ) - \Rot^{\bm{T}}\fEr(t) \right] = \bm{0}, \hspace{1cm} \forall \ \delb\qrr \in \R^m.
 \end{equation}
Using Eq.\eqref{eq:redForcTens}, we can retrieve the transformed ROM equations and transformation laws for internal force tensors, reduced mass and damping matrices:
\begin{equation}
    \mathbf{\tilde{\tilde{M}}} \qrrdd + \mathbf{\tilde{\tilde{C}}}  \qrrd + \Korr \cdot \qrr + \Ktwrr \cddot (\qrr \otimes \qrr) + \Kthrr \cdddot\   (\qrr \otimes \qrr \otimes \qrr ) = \Rot^{\bm{T}}\fEr(t)
\end{equation}
 with
 \begin{subequations}
\begin{align}
    \Korr &= \RotT \Kor \Rot \\ 
    \Ktwrr &= \left( \RotT \Ktwr \cdot \Rot \right)\tensCon{2}{1} \\ 
    \Kthrr &= \left( \left( \RotT \Kthr \cdot \Rot \right) \tensCon{3}{1} \Rot \right) \tensCon{2}{1} \Rot\\
    \mathbf{\tilde{\tilde{M}}} &= \RotT \Mr \Rot  \\
    \mathbf{\tilde{\tilde{C}}} &= \RotT \Cr \Rot.
\end{align}
\label{eq:transfTensors}
\end{subequations}
If we choose $\Rot$ such that
\begin{equation}
\V\Rot = \W,
\label{eq:transfUW}
\end{equation}
the ROM in Eq.\eqref{eq:transfTensors} is written for the RB $\W$ since
\begin{equation}
    \q \approx \V\qr = \V\left(\Rot \qrr \right) = \W \qrr.
\end{equation}
Given $\W$ and $\V$, the coefficients of the linear transformation stored in matrix $\Rot$ are uniquely defined (since $\W$ and $\V$ must span the same subspace) and can be computed by projecting Eq.\eqref{eq:transfUW} onto $\V$, solving for $\Rot$:
\begin{equation}
    \Rot = \left(\VT \V\right)^{-1}\VT \W.
    \label{eq:transfMtrC}
\end{equation}
With $\Rot$ available, Eqs.\eqref{eq:transfTensors} are used to retrieve the transformed tensors written for the orthogonalized RB $\W$.
The reader should notice that transformation formulas are still valid when applied to tensors in the sparse form obtained from EED identification.
However, application of this transformation to sparse tensor does not produce transformed sparse tensors. 
To regain the advantages of a sparse representation of tensors, terms multiplying the same monomials in the polynomial expression of the reduced forces can be collected together and saved in a sparse tensor format more suitable for ROM time integration.
%sparsification algorithm presented in Appendix \tcg{scrivi alg.} can be applied to the transformed tensors, before the ROM is used for time integration.\\ 
%The procedure described in this section for tensor identification using EED-ECSW is summarized in Algorithm \ref{alg:EEDnew}.
%algorithm for reduced mesh
\begin{algorithm}[H]
\caption{EED-ECSW for ROM Tensor Identification}
\label{alg:EEDnew}
\textbf{Input:} \texttt{FeModel}; RB $\V$; ECSW reduced element set $\Elsr$; ECSW element weights $\xielvec^*$; linear stiffness of elements in $\Elsr$ $\{\mathbf{K_{1}},...,\mathbf{K_{\nelr}}\}$ 
\\
\textbf{Output:} ROM tensors $\Korr,\Ktwrr,\Kthrr$ 
 \begin{algorithmic}[1] 
    \Statex \textit{RB orthogonalization}
    \State $\W \gets \texttt{orth}(\V)$\sscr{a}
    \Statex \textit{Tensor identification with EED}
    \State $ \{\qj{1},...,\mathbf{q}^{(N_q)}\} \gets 
    \texttt{displEED($\V$)}$\sscr{b}
    \State Assemble reduced FE model \texttt{FeRed}, restricting \texttt{FeModel} to elements in $\Elsr$
    \State $\mathbf{S_{\tilde{K}^{t}}} \gets \texttt{tuple}(N_q)$\sscr{c}
    \For{$i \in \{1,...,N_q\}$}
    \State  $\{\mathbf{K_1^t},...,\mathbf{K_{\nelr}^t}\} \gets \texttt{FeRed.ElementStiffNlin}(\qj{i})$\sscr{d}
    \State $\Ktgrnl \gets \texttt{ECSWmodel.tangStiff($\xielvec^*,\{\mathbf{K_1^t},...,\mathbf{K_{\nelr}^t}\}$},\{\mathbf{K_{1}},...,\mathbf{K_{\nelr}}\},\V)$\sscr{e}
    \State $\mathbf{S_{\tilde{K}^{t}}}\{i\} \gets \Ktgrnl$
    \EndFor
    \State $[\Ktwr,\Kthr] \gets \texttt{tensorID($\mathbf{S_{\tilde{K}^{t}}},\V$)}$ \sscr{f}
    \Statex \textit{Transform tensors from RB $\V$ to RB $\W$}
    \State $\Rot \gets \left(\VT \V\right)^{-1}\VT \W $ \sscr{g}
    \State $ [\Korr,\Ktwrr,\Kthrr] \gets \texttt{transformTensors($\Kor,\Ktwr,\Kthr,\Rot$)}$\sscr{h}
    \State $ [\Korr,\Ktwrr,\Kthrr] \gets \texttt{sparsifyTensors($\Korr,\Ktwrr,\Kthrr$)}$\sscr{i}
\end{algorithmic}
\end{algorithm}
\footnotesize
\noindent \sscr{a}\texttt{orth}($\V$) function orthogonalizes RB $\V$.\\
\sscr{b}\texttt{displEED}$(\V)$ returns a set of $N_q$ nodal displacements for tensor identification with EED.  $N_q = (m^2+5m)/2$, where $m$ is the \indent number of column vectors in $\V$.\\
\sscr{c}\texttt{tuple}$(i)$ function initializes a tuple array with $i$ elements.\\
\sscr{d}\texttt{FeRed.ElementStiffNlin($\qj{i}$)} returns the nonlinear tangent stiffenss matrix for all elements in the \texttt{FeRed} model, for \indent imposed displacement $\qj{i}$.\\
\sscr{e}\texttt{ECSWmodel.tangStiff} computes the approximation of the reduced tangent stiffness with ECSW, as in Eq.\eqref{eq:approxECSW}.\\
\sscr{f}\texttt{tensorID} function solves Eqs.\eqref{eq:identifEqs},\eqref{eq:identifEqs1} for the unknown tensor coefficients.\\
\sscr{g}Transformation matrix $\Rot$ is computed as in Eq.\eqref{eq:transfMtrC}.\\
\sscr{h}\texttt{transformTensors} function transforms the tensors implementing the operations in Eq.\eqref{eq:transfTensors}.\\
\sscr{i}\texttt{sparsifyTensors} function converts tensors from a dense to a sparse representation.\\
\normalsize

%case study
\section{Applications} \label{sec:applications}
The proposed methodology has been validated on two different structures subjected to random loads: a curved panel and a nine-bay aeronautical panel.
We discuss these examples in two dedicated sections, drawing a comparison between the performances of EED-ECSW and standard EED identification.
%Implementation details
To demonstrate the non-intrusive capabilities of the presented method, we construct the ROM based on ABAQUS commercial FE package.
MATLAB programming language was used for ROM construction and ROM time integration.
Practical implementation details for ROM construction with ABAQUS are illustrated alongside the presentation of the curved plate, in the following section \ref{sec:curvedPanel}.
\subsection{Curved Panel}
\label{sec:curvedPanel}
\begin{figure}
    \centering
\includegraphics[width=0.9\linewidth]{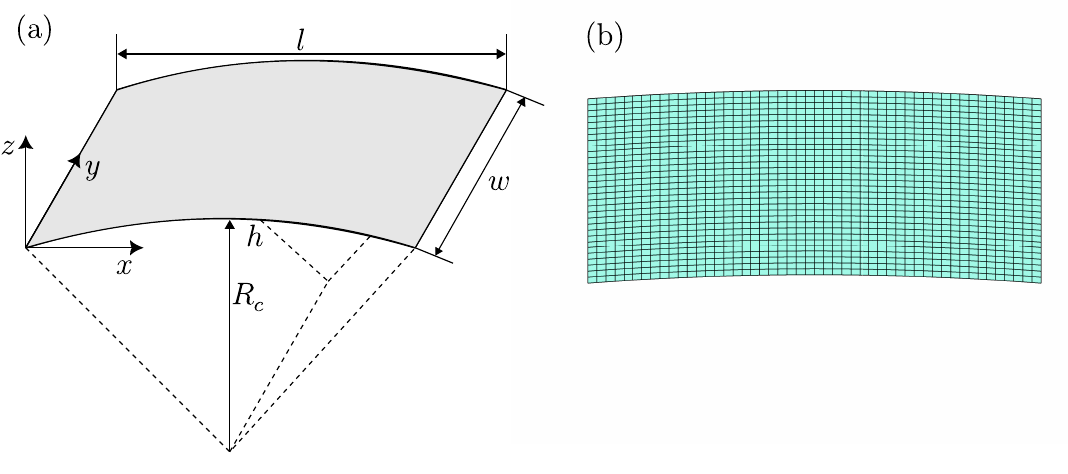}
    \caption{Geometry (a) and FE mesh (b) of the curved panel.}
    \label{fig:curvPanGeoMesh}
\end{figure}

\begin{figure}
    \centering
    \includegraphics[width=1\linewidth]{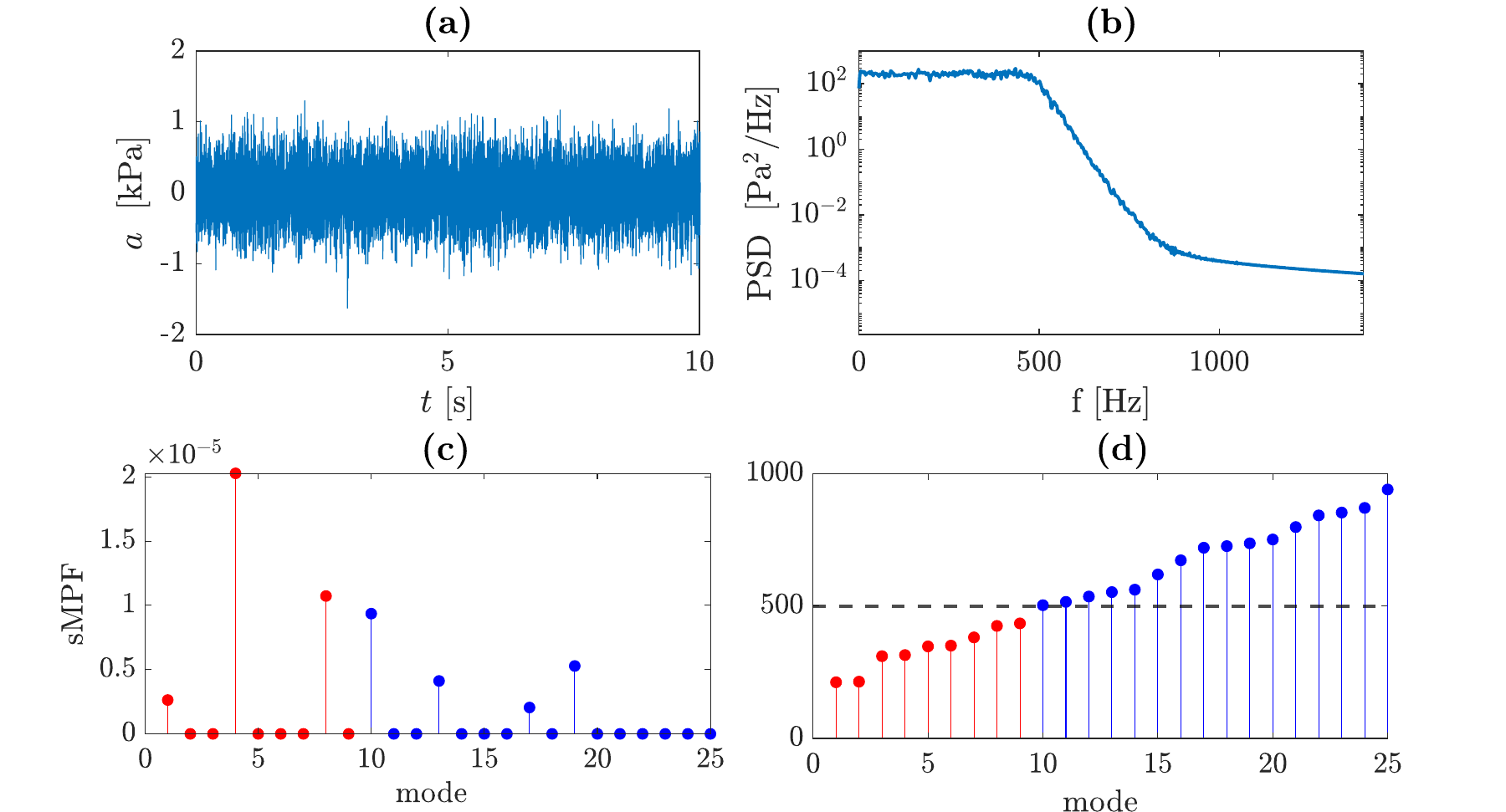}
    \caption{Time history (a) and PSD (b) of the uniform in space pressure applied to the curved panel.
    PSD was obtained with Welch's method.
    In (c) static modal participation factor of the pressure load, while in (d) natural frequencies of the curved panel. VMs before $f_{\text{coff}}$ are plotted in red, the one after in blue.}
    \label{fig:curvPanPress}
\end{figure}

As a first example we present a slightly curved rectangular panel of length $l = 0.4\  \text{m} $, width $w = 0.25\ \text{m}$ and thickness $\text{t} = 0.8 \ \text{mm}$.  
The height of the middle line in the length direction is $h = 0.0079 
\ \text{m}$, which corresponds to a curvature radius $R_c =  2.54\  \text{m}$.
The panel is clamped at all edges.
A sketch of the geometry is provided in Fig. \ref{fig:curvPanGeoMesh}a.
The material is assumed to be linear elastic with elasticity modulus $\text{E} = 70 \ \text{GPa}$, Poisson's ratio $\nu \ = \ 0.3$ and density $\rho = 2700 \ \text{Kg}/\text{m}^3$.
The FE model was constructed with $1550$ S4R elements provided by ABAQUS (mesh in Fig. \ref{fig:curvPanGeoMesh}b), resulting in $8820$ unconstrained degrees of freedom.
\\
A uniform in space pressure, varying randomly in time as a white noise band limited process in the frequency range $0-500\ \text{Hz}$, with an \textit{Overall Sound Pressure Level} (OASPL) of $144 \ \text{dB}$ ($P_{ref}$ = 20 $\mu\text{Pa}$), is applied to the top surface of the panel. 
This load aims to mimic the effect of acoustic loading on aerospace structures \cite{Hollkamp2018}.
The time varying load applied to the FE model was constructed by multiplying the nodal load equivalent to a uniform unitary pressure $\mathbf{p}$, by a time varying amplitude function $a(t)$ as
\begin{equation}
    \fe(t) = a(t) \mathbf{p}.
    \label{eq:feForcForm}
\end{equation}
The amplitude function was constructed with samples from the standard normal Gaussian distribution, filtered with  a low-pass Butterworth filter (order $12$, cut off frequency $f_{\text{coff}}$ = $500$ Hz), and eventually re-scaled to match the target OASPL, adopting the same approach as in \cite{Schoneman2017}.
As can be seen from Fig. \ref{fig:curvPanPress}a-b, the estimated PSD of the constructed pressure signal is flat up to $f_{\text{coff}}$ and rapidly drops thereafter, as desired.\\
Rayleigh damping was adopted, with values of $\alpha = 32.96$ and $\beta = 1.162 \cdot 10^{-5}$ obtained by imposing a modal damping ratio of $0.02$ to the first 5 VMs (in a least square sense). 
\subsubsection{Reduced order basis construction}
Vibration Modes computation is the first step for the construction of the RB. 
The eigenvalue problem in Eq.\eqref{eq:eigLin} was solved with ABAQUS, extracting the first $25$ VMs and associated natural frequencies.
Vibration modes selection for forced response was performed based on a \textit{Static Modal Participation Factor} (SMPF) defined as 
\begin{equation}
    \text{SMPF}_i = \frac{\bm{\phi}_i^T \mathbf{p}}{\bm{\phi}^T_i \Ko \bm{\phi}_i}\|\bm{\phi}_i\|_2.
    \label{eq:MPF}
\end{equation}
As shown in Appendix A, $\text{SMPF}_i$ provides the contribution of VM $i$, to the linear static solution to load vector $\mathbf{p}$, for  unit normalized modes.
As opposed to classical formulations for the modal participation factor, this measure takes into account the increasing in modal stiffness with the mode number and thus penalizes stiffer modes.
The employment of the sMPF is justified by the fact that the spectogram of the amplitude function is flat over broad-band and thus a ranking of the VMs can only be made based on their shapes. 
As such, by plotting the sMPF for the computed modes (Fig. \ref{fig:curvPanPress} c,d), we decided to include in the RB VMs $1,4,8,10,13$,$17$ and 19.
These VMs have non zero MPF and their associated natural frequencies fall either within (VMs $1,4,8$) or just outside (VMs $10,13,17$ and $19$) the exicitation bandwidth. Therefore they are expected to contribute to the dynamic response, as a-posteriori shown in section \ref{sec:curvedPanConv}.\\
With the VMs available, SMDs were constructed by multiple solution of the linear problem in Eq.\eqref{eq:SMDdef}.
To this purpose, the directional derivatives of the tangent stiffness matrix were computed with central finite differences as in Eq.\eqref{eq:SMDsRHS} from the Abaqus outputs of the tangent stiffness in the direction of the VMs.
Inclusion of active VMs ($7$ vectors) and all the corresponding SMDs ($28$ vectors) led to a RB of size $35$.
\subsubsection{ECSW model construction}
\begin{figure}
    \centering
    \includegraphics[width=0.7\linewidth]{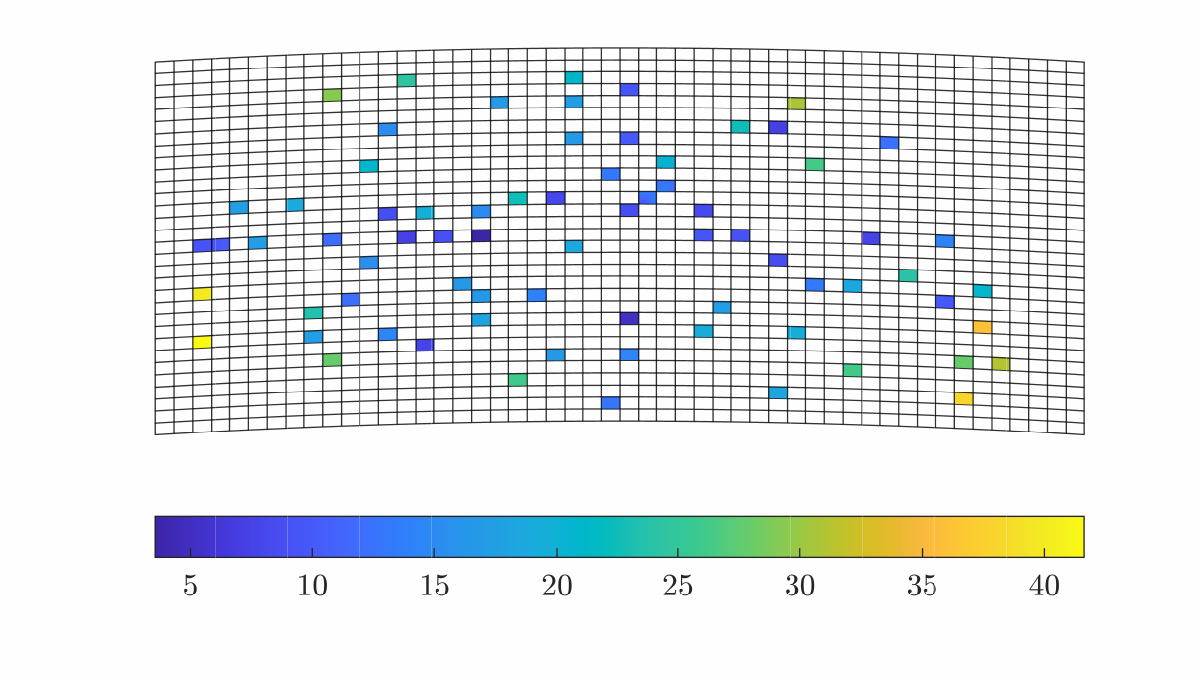}
    \caption{Reduced mesh for the curved panel. Color intensity maps to ECSW element's weight. The number of active elements is $\nelr = 73$, corresponding to $4.71\%$ of the total number of elements in the original FE mesh.}
    \label{fig:redMeshCurved}
\end{figure}

The Energy Conserving Sampling and Weighting hyperreduced model was computed using the scheme presented in section \ref{sec:ECSWmodel} and summarized in algorithm \ref{alg:redMesh}.
As a first step, the SQM for ECSW training was constructed with all the VMs in the RB. 
Latin Hypercube Sampling was used to generate $\Nt= 50$ quasi random samples of the VMs amplitude vector $\bm{\gamma}$, with $N_t = 45$ of them used for training and $N_v = 5$ for validation.
Limits of the Latin Hypercube were enforced by setting the bounding parameter $\alpha$, defined in section \ref{sec:SQM_4train}, to $\alpha = 0.6\  \text{t}$, where  $\text{t}$ is the thickness of the structure.
This choice was driven by experience and by comparison between linear and nonlinear internal forces for various imposed displacement levels, making sure that nonlinearities are activated.
Abaqus was used to compute the element level nonlinear internal forces for the displacements generated with the SQM and to extract the linear elemental tangent stiffness matrices.
Matrices $\Gecsw$ and $\becsw$ are assembled from these forces after subtraction of their linear part, and by projection on the physical RB $\V$ (see Eq, \eqref{eq:forceSetNl} and \eqref{eq:ECSWGb}).
With these quantities available, the NNLSQ algorithm was employed with a relative tolerance $\tau = \ 1 \cdot 10^{-3}$ to solve the optimization problem in Eq.\eqref{eq:minProbl}, retrieving the ECSW reduced mesh and associated weights $\bm{\xi}^*$.
The computed reduced mesh, displayed in Fig. \ref{fig:redMeshCurved}, featured $73$ elements (corresponding to $4.71 \%$ of the total number of elements) and a validation error of $\epsilon_{ECSW} = 5.6 \cdot 10^{-4}$, which we deemed acceptable for identification.
\subsubsection{Tensor identification with EED-ECSW}
\label{sec:tensIdcurv}
Tensor identification was performed according to the EED scheme, which
requires the evaluation of the reduced nonlinear tangent stiffness matrix for $665$ different displacements.
For each enforced displacement, the nonlinear tangent stiffness was computed only for the elements in the reduced mesh.
This was accomplished by creating an Abaqus input file that includes only the elements in the reduced mesh and their corresponding nodal imposed displacements, then submitting it for analysis.
The differences between nonlinear and linear reduced tangent stiffness were then approximated using the ECSW weights as in Eq.\eqref{eq:approxECSW}, in which the physical RB $\V$ was used for projection.
Notice that this formula requires the knowledge of linear element stiffness, previously computed in the identification of the reduced mesh.
Eventually, tensor entries have been identified by solving a set of algebraic linear systems of equations obtained from Eq.\eqref{eq:identifEqs} and \eqref{eq:identifEqs1}, written with the approximated tangent stiffness.
\par
The tensors identified in the physical RB $\V$ were eventually  transformed to their counterpart corresponding to a mass orthormalized RB $\W$, obtained from $\V$ using the Gram-Schmidt algorithm, following the procedure presented in section \ref{sec:tensRot}.
\subsubsection{Construction time}
\label{sec:constrTimeCurvedPanel}
In order to quantify the advantages of the proposed approach over standard EED, we also constructed the tensors with this second method.
A comparison of the computational times for these two approaches applied to the curved panel ROM are reported
in table \ref{tab:compTimeCurvP}. 
Tensor construction has been performed on the Euler cluster of ETH Z\"urich, using 10 cores, with 3 Gb RAM each.
Computational times for RB construction and tensor identification are reported in table \ref{tab:compTimeCurvP}.    \\
By using EED-ECSW instead of standard EED for tensor identification, the overall ROM construction time decreased from $1908.5$ s to $599.9 $ s, corresponding to a speed up of 3.18.
The overall computational advantage is modest in this case due to the small dimensions of the underlying FE model. \par
In order to understand which operations are accelerated by using the reduced mesh,  
we tracked the computational time for different processes involved in tensor identification. 
Specifically, we categorized the computational times into six different categories, reported under the keywords 'Reduced Mesh', 'Abaqus', 'Reading', 'Identification', 'Tensor Transformation' and 'other'.\\ 
The computation of the reduced mesh is an overhead cost that is sustained only in EED-ECSW.
The cost for the reduced mesh is dominated by the computation of the training forces with the FE model (32.0 s), while solving the optimization problem in Eq.\eqref{eq:minProbl} with NNLSQ comes almost with no additional cost (0.14 s).\\
Besides the computation of the reduced mesh, all the other operations are performed both within EED and EED-ECSW. 
The computational time reported under the  keyword 'Abaqus' comprises the time spent by Abaqus to run the analysis and return the tangent stiffness for the queried displacements.
Since the reduced mesh features less elements than the FE mesh, the computation cost for EED-ECSW is lower than the one for EED. Notice however, that the speed up of $3.44$ for these operations is not equal to the ratio between the number of elements in the mesh and in the reduced mesh ($\nel/\nelr = 21.23$), as one would expect.
This is because the model is relatively small, and the cost of some sub-operations in the Abaqus analysis does not scale proportionally with the number of elements (e.g. opening and closing files for writing the matrices, communication with OS), but is rather constant.\\
Under the keyword 'Reading' we report the computational time spent by the Matlab code to read the tangent stiffness matrices saved by Abaqus in .mtx files.
Since EED-ECSW uses the reduced mesh for tensor identification, the number of entries of the tangent stiffness is limited as compared to EED, and reading is thus faster, even when sparsity of the full FE tangent stiffness matrix is exploited.
Specifically, for this example the recorded speed up in reading operations was $4.97$ between ECSW-EED and EED. \\
Small differences in the two methods are observed for the computation times reported as 'Identification' and 'Tensor transformation'. 
Under the first keyword we include the time spent for solving the equations to identify tensors coefficients, while under the second keyword the time spent to transform the tensors from the physical RB to the orthogonalized RB. 
In fact, the computation cost for these operations does not scale with the number of elements in the FE formulation yet only with the size of the ROM, resulting in marginal differences between the two methods.\\
Lastly, we report under the keyword 'other', the computational times for the operations not included in the previous fields, such as array initializations, Abaqus input file writing and file saving.
Although the cost for these operations are small, a speed up of $5.3$ is observed with EED-ECSW.
\begin{table}[h!]
\centering
\vspace{0.2cm}
\textbf{Reduction basis construction}\\
\vspace{0.1cm}
\begin{tabularx}{0.8\textwidth}{>{\raggedright\arraybackslash}X 
   >{\raggedright\arraybackslash}X 
   >{\raggedright\arraybackslash}X
   >{\raggedright\arraybackslash}X }
   \hline
    VMs & $12.3$ s &\   & \ \\
    %\hline
    SMDs & $36.1$ s &\   &\  \\
    \hline
    \textbf{total RB} & $\bm{48.4}$\ \textbf{s} & \  &\  \\
\end{tabularx}\\
\vspace{0.2cm}
\textbf{Tensor identification}\\
\vspace{0.1cm}
\begin{tabularx}{0.8\textwidth}{>{\raggedright\arraybackslash}X 
   >{\raggedright\arraybackslash}X 
   >{\raggedright\arraybackslash}X
   >{\raggedright\arraybackslash}X }
  \hline
    &EED  & EED-ECSW  & speed-up\\
    \hline
    reduced mesh & - & 32.1 s & -\\
    %\hline
    Abaqus & 1458.8 s & 424.4 s &3.44 \\
    %\hline
    reading & 376.7 s & 75.76 s &4.97 \\
    %\hline
    identification & 17.6 s & 17.8 s &0.99 \\
    %\hline
    transformation & 0.13 s & 0.11 s &1.18 \\
    other &        6.9 s & 1.3 s & 5.30\\  
    %\hline
    \hline
    \textbf{total Id.} & $\bm{1860.1}$\ \textbf{s} & $\bm{551.47}$\ \textbf{s} &$\bm{3.37}$\\
    \hline  \hline
    \textbf{total RB \& Id.} & $\bm{1908.5}$\ \textbf{s} & $\bm{599.9}$\ \textbf{s} &$\bm{3.18}$\\
    \hline
\end{tabularx}
\caption{Tensors construction times for the curved panel ROM: comparison between EED and EED-ECSW.}
\label{tab:compTimeCurvP}
\end{table}
\subsubsection{Time integration and results}
\label{sec:curvedPanRes}
To evaluate accuracy of the proposed ROM, we integrated in time the ROM with tensors identified by standard EED, the ROM with tensors identified by EED-ECSW, and finally, the HFM.
We performed time integration for the ROMs in Matlab, using a Matlab subroutine that implements the Newmark-$\beta$ integration scheme \cite{geradin2015mechanical}, while time integration of the HFM was performed in Abaqus using the HHT scheme.
In both cases, the time step was set to $4.167 \cdot 10^{-5}\ \text{s}$ and the analysis was run for $10\ \text{s}$ of real time simulation. 
Time histories of nodal degrees of freedom were recorded for two nodes identified in the following with 'node A' and 'node B', with respective coordinates $(x_A = 0.5\cdot l  , y_A = 0.516\cdot w, z_A = h)$ and  $(x_B = 0.34\cdot l, y_B = 0.322 \cdot w, z_B =0.90 \cdot h)$. Each time domain signal had 240,002 data points.
%Displacements time series obtained with the Abaqus HFM, with ROMs constructed with EED and EED-ECSW and eventually with a ROM for the linearized system (with 30 VMs in the RB) are  shown in Fig. \ref{fig:displCurved}.
%From the figure it can be seen that both ROMs are in good agreement with the HFM, and that  displacement level of the nonlinear vibration is generally larger than the linear.
The Power Spectral Density for displacements degrees of freedom of these nodes was estimated using the Welch's method \cite{hayes1996statistical} with segments of 9,000 sample points and $50\%$ overlap, resulting in a frequency resolution of $2.67$ Hz. \par
\begin{figure}
    \centering
    \includegraphics[width=0.7\linewidth]{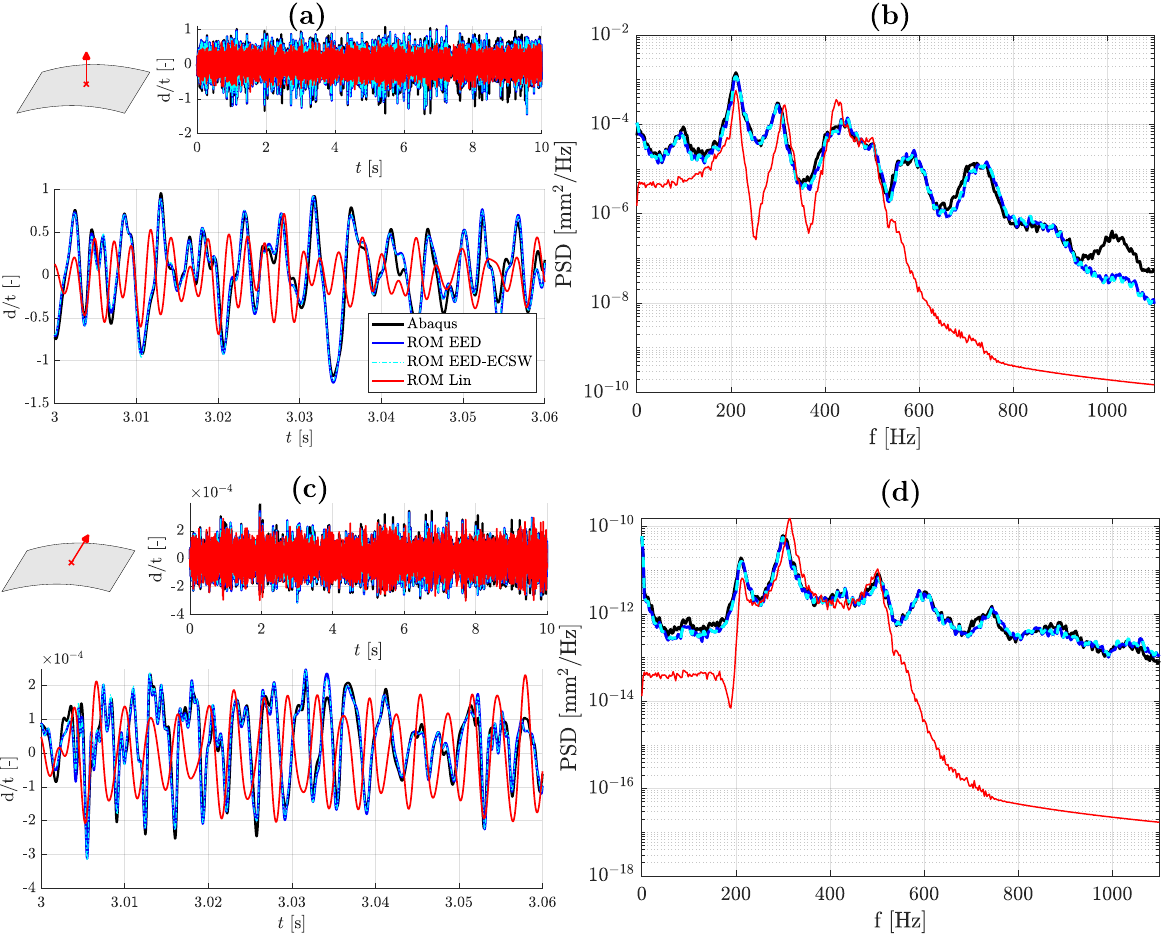}
    \caption{Time domain displacements of node A $(x_A = 0.5\cdot l  , y_A = 0.516\cdot w, z_A = h)$ and their PSDs. In (a) and (b) out of plane displacements, while in (c) and (d) in plane displacements.
    Displacements in (a) and (c) are normalized with respect to the thickness of the plate.
    The four lines in each plot correspond to solutions obtained using Abaqus HFM (Abaqus), the ROM with tensors identified through EED (ROM EED), the ROM with tensors identified using EED-ECSW (ROM EED-ECSW) and the linear ROM (ROM Lin).  The RB used in the nonlinear ROMs consists of 7 VMs and 28 SMDs.}
    \label{fig:displCurved1}
\end{figure}
\begin{figure}
    \centering
    \includegraphics[width=0.7\linewidth]{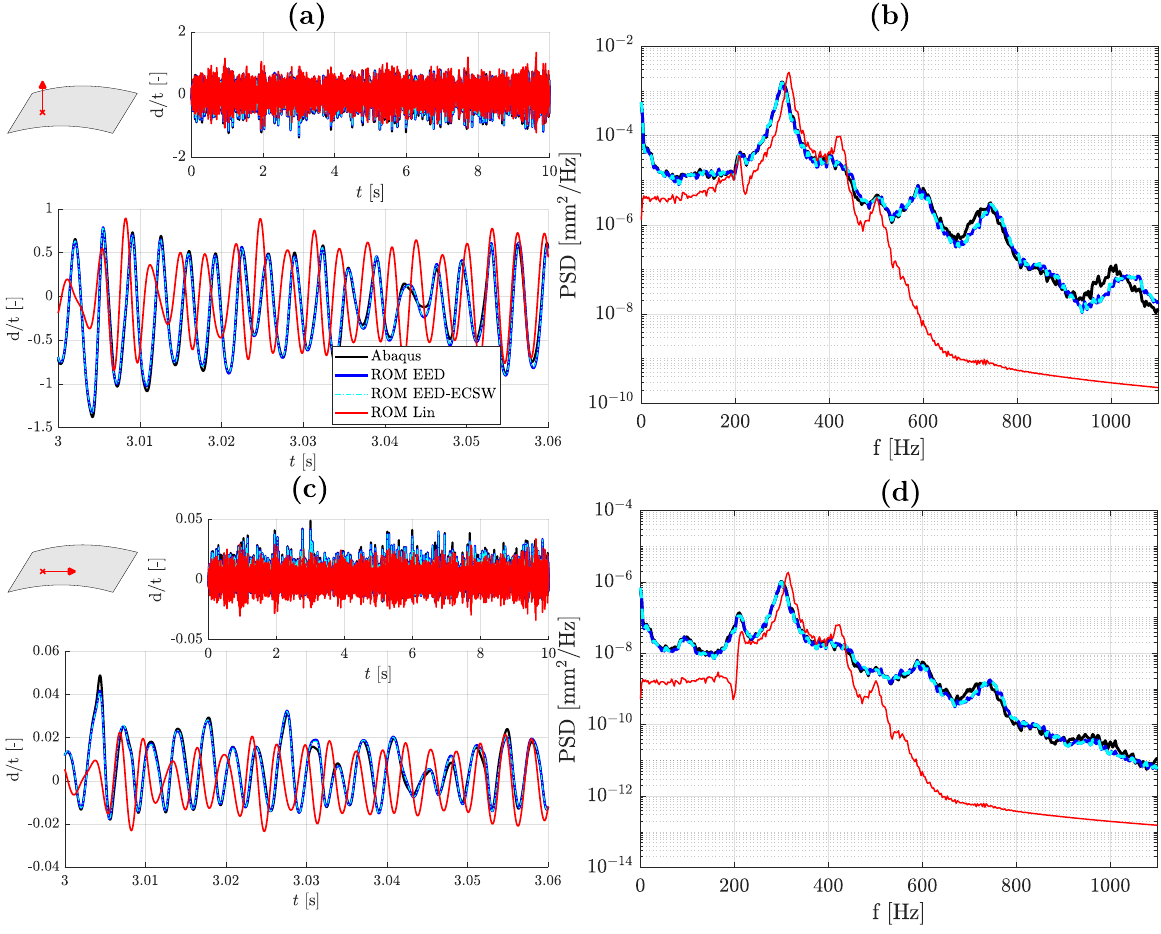}
    \caption{Time domain displacements of node B $(x_B = 0.34\cdot l, y_B = 0.322 \cdot w, z_B =0.90 \cdot h)$ and their PSDs. In (a) and (b) out of plane displacements, while in (c) and (d) in plane displacements.
    Displacements in (a) and (c) are normalized with respect to the thickness of the plate.
    The four lines in each plot correspond to solutions obtained using Abaqus HFM (Abaqus), the ROM with tensors identified through EED (ROM EED), the ROM with tensors identified using EED-ECSW (ROM EED-ECSW) and the linear ROM (ROM Lin). The RB used in the nonlinear ROMs consists of 7 VMs and 28 SMDs.}
    \label{fig:displCurved2}
\end{figure}

In Fig. \ref{fig:displCurved1} and Fig. \ref{fig:displCurved2} we plot the time histories and PSDs of displacements of node A and node B respectively.
The ROMs are in excellent agreement with the Abaqus HFM, both in the time and frequency domain. 
A slight mismatch is observed after long simulation time in the time histories and for frequency larger than 900 Hz in the PSD (see Fig. \ref{fig:displCurved1} b).
For comparison purposes, we plotted the solution to the same load obtained using a ROM constructed for the linearized FE model (with 30 VMs in the RB) alongside the response of the nonlinear models.
By comparing the PSD of the nonlinear solution with the linearized solution, it is easy to notice that geometric nonlinearities introduce a strong static component in the response, a smearing of the PSD peaks, as well as a pronounced excitation of frequencies outside the load frequency bandwidth.
Moreover, PSD peaks are shifted to lower frequencies, suggesting a softening behavior of the structure, as is well known for slightly curved panels.\\
The total recorded computational time for integrating the Abaqus HFM was $73.22$ h on the Euler cluster of ETH Z\"urich, using 10 cores, with 3 Gb RAM each. 
The ROM with EED-ECSW tensors was integrated on the same machine in only $655.3$ s, achieving an online speed up factor of $402$. 
A similar time integration speed up was obtained for the ROM with tensors constructed using EED, as expected.
\subsubsection{Convergence analysis for number of VMs in RB}
\begin{figure}
    \centering
    \includegraphics[width=1\linewidth]{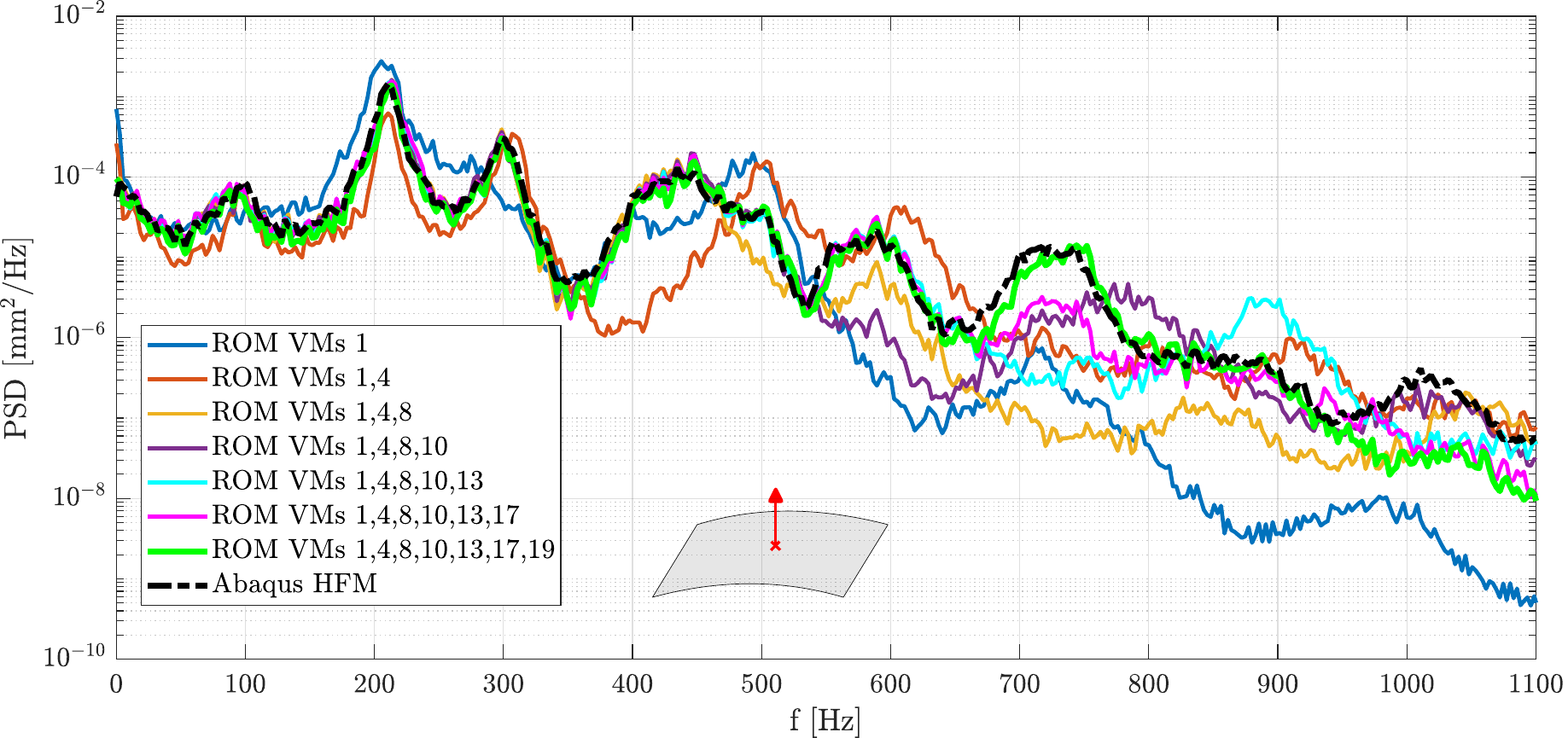}
    \caption{Convergence analysis for the curved panel by increasing the number of VMs in the ROM. PSD of nodal displacement along $z$ axis of node A.}
    \label{fig:convAnalysis}
\end{figure}

In this subsection we present a convergence analysis of the PSD obtained using the ROM, for an increase in the number of VMs retained in the RB. \\
Eight different ROMs were computed for different RBs, subsequently constructed by including VMs with non-zero sMPF, for increasing eigenfrequency.
Specifically, we progressively included in the RBs VMs $1,4,8,10,13,17$ and 19, along with their corresponding SMDs, one at a time.\\
Convergence of the PSD of displacement history of node A along $z$ is displayed in Fig. \ref{fig:convAnalysis}. 
It is worth to notice that multiple VMs are needed for a correct representation of the PSD over all the frequency span.
In particular, addition of single VMs improves the PSD results in the neighbourhood of the mode eigenfrequency.\\
The requirement of multiple VMs in the RB for correct PSD estimation leads to large ROMs whose computation is generally intensive, hence justifying the need for the efficient ROM identification scheme herein presented.
\label{sec:curvedPanConv}
%Nine Bay
\subsection{Nine-Bay Panel}
As a second test case we present the application of our method to a model of the aircraft sidewall fuselage panel experimentally investigated in \cite{buehrle2000finite}. 
The nine bay panel consists of a skin panel reinforced with a frame and a longeron substructure dividing the skin into nine bays.
The response of this panel to acoustic loading was assessed in \cite{przekop2006nonlinear} and \cite{Perez2014} using projection ROMs.
In the first study \cite{przekop2006nonlinear}, the RB was constructed using only VMs, with tensors identified through ED. 
In contrast, the second study \cite{Perez2014} included Dual Modes in the RB, and EED was used for tensor identification.
The FE model that we used is derived from the one employed in \cite{przekop2006nonlinear} by refining the mesh, and by substituting beam elements connecting the shell surfaces with strings of shell elements.
Mesh refining was performed to demonstrate the applicability of the proposed approach to large FE models, while the substitution of beam elements with shells was intended to create a FE model composed of elements of the same type, simplifying the implementation.
However, the proposed method for ROM construction could also be applicable to FE models featuring different element types. \\ \par
The dimensions of the skin rectangular panel are $1.4760\  \text{m}$ by $0.6365\ \text{m}$, while 
the thickness for all shell surfaces is $1.3\ \text{mm}$.
The adopted linear elastic material has Young's modulus $E = 72.395 \ \text{GPa}$, Poisson's ratio $\nu = 0.33$ and density $\rho = 2794  \ \text{Kg/m}^3$.
The FE model was constructed with Abaqus, using $33,208$ S4R elements and features $200,886$ unconstrained degrees of freedom.
A sketch of the FE mesh is provided in figure \ref{fig:meshNineBay}.a.
The nine-bay panel is constrained by blocking all the displacements degrees of freedom at the edges of the skin, leaving the nodal rotations free.
The applied load is a uniform pressure in space, varying randomly in time as a band limited white noise, acting on the skin of the panel.
The load was constructed as described in section \ref{sec:curvedPanel} for the curved panel, imposing an OASPL of $147$ dB and a cut-off frequency of $500$ Hz.\\
Rayleigh damping model was adopted, with values of $\alpha$ and $\beta$ determined by imposing a damping ratio of $0.01$ to the first two structural VMs.
\subsubsection{Reduced order model}
\begin{figure}
    \centering
    \includegraphics[width=0.9\linewidth]{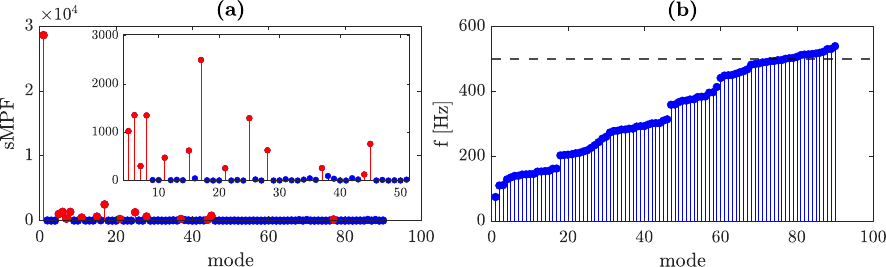}
    \caption{sMPF (a) and natural frequencies (b) of the nine-bay panel. In (a), the 15 VMs included in the RB are plotted in red. Seventy-six  VMs have frequency less than the $f_{\text{coff}}$ (dashed line in (b)).}
    \label{fig:sMPFNineBay}
\end{figure}

The first $90$ VMs and natural frequencies were computed with Abaqus.
Modal selection for this problem is crucial, since $76$ VMs fall within the frequency bandwidth up to the cut off frequency of the load (Fig.\ref{fig:sMPFNineBay}.b), and their blind inclusion in the RB would result in an inefficient ROM.
As such, the RB was constructed from the 15 VMs with the highest sMPF (see Fig.\ref{fig:sMPFNineBay}.a).
Subsequently, SMDs were computed non-intrusively from these VMs, and ranked based on a newly proposed criteria that builds on the same spirit of the selection strategy proposed in \cite{tiso2011optimal}.
In particular, the SMPFs of the VMs included in the RBs were cross multiplied, obtaining
\begin{equation}
    {\text{MDPF}}_{ij} = \text{SMPF}_{i} \cdot \text{SMPF}_{j} \quad \text{for}\ \  i \geq j,
\end{equation}
where $i,j$ are indices corresponding to the VMs in the RBs, and $\text{MDPF}_{ij}$ is name the \textit{Modal Derivative Participation Factor} of the $i$\textsuperscript{th} and $j$\textsuperscript{th} VMs.
The rationale behind the design of this ranking measure is that the sMPF for a given VM is assumed to be proportional to the modal displacements of that VM when the load is applied to the linearized system. 
If this holds, this ranking criteria for SMDs is similar to the selection criteria proposed in \cite{tiso2011optimal}.\\
Based on this, we form the RB by adding to the  selected $15$ VMs the $35$ SMDs with the highest values of MDPF index.
\\  \par
\begin{figure}
    \centering
    \includegraphics[width=1\linewidth]{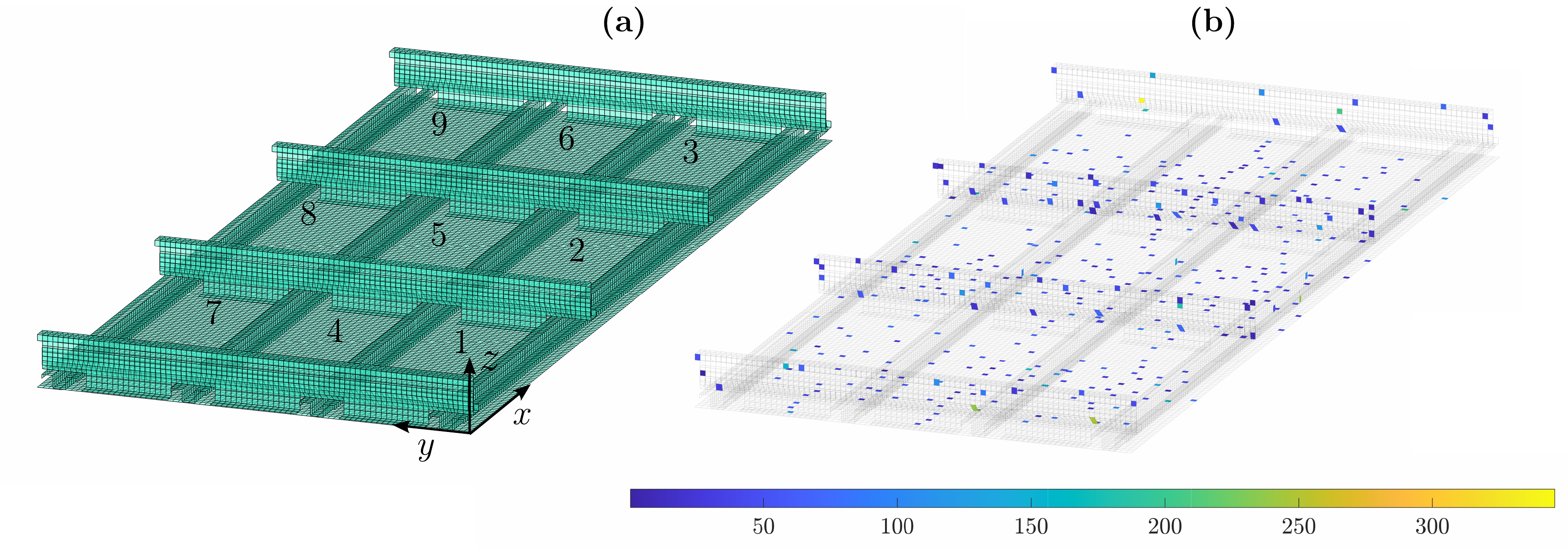}
    \caption{In (a) FE mesh for the nine-bay panel: the bays on the skin are numbered from 1 to 9. In (b) its reduced mesh computed for $\tau = 0.001$. Intensity color is proportional to element weights. The hyperreduced mesh consists of only 395 elements, corresponding to 1.19 $\%$ of the elements in the FE mesh.}
    \label{fig:meshNineBay}
\end{figure}

The constructed RB consisted of $50$ vectors requiring $1325$ evaluations of the tangent stiffness matrix within the EED tensor identification procedure.
Reduced order model tensors were identified with EED-ECSW and, for comparison purposes, with standard EED.
The SQM displacements, used for reduced mesh construction in EED-ECSW identification, were obtained from $\Nt=70$ quasi random samples of the VMs amplitude vector (only VMs in RB), $60$ of which were used for training and the rest for validation. 
Limits of the Latin Hypercube were enforced by setting the bounding parameter $\alpha$ (see section \ref{sec:redMeshComp}) equal to the thickness of the structure.
The optimization problem for ECSW mesh computation was solved for relative tolerance $\tau = 0.001$, retrieving a reduced mesh with $395$ elements (corresponding to $1.1895\%$ of the total number of elements in the FE model) and with validation error $ \epsilon_{\text{ECSW}} = 0.0024$.
The reduced mesh is displayed in Fig. \ref{fig:meshNineBay}.b.
 Tensors were identified with EED-ECSW with the reduced mesh available, following the same procedure adopted in section \ref{sec:tensIdcurv} for the curved panel.

\begin{figure}
    \centering
    \includegraphics[width=1\linewidth]{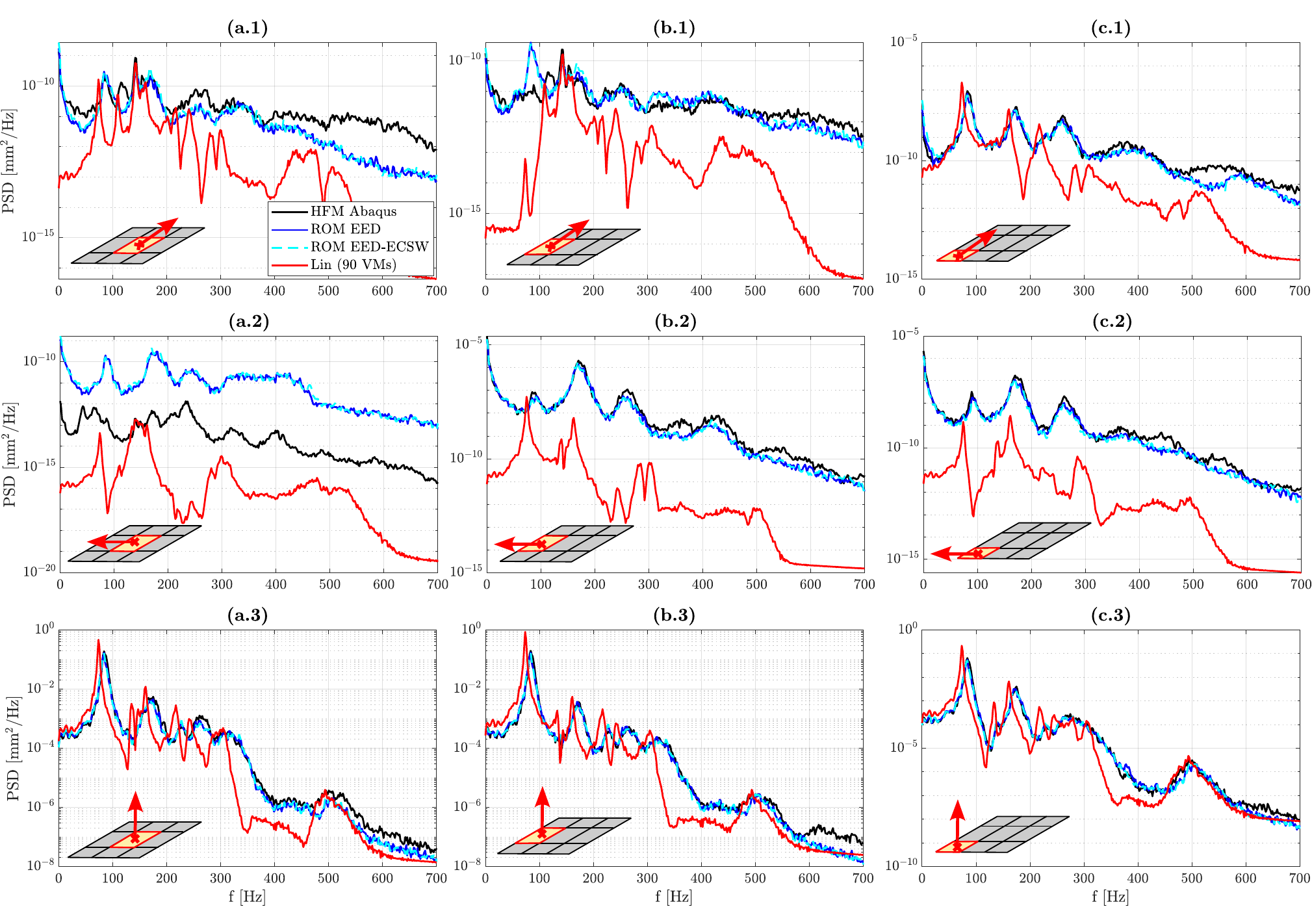}
    \caption{Power spectral density of displacements for three different nodes located in the middle of bay 5 (a.1-a.3), of bay 8 (b.1-b.3) and of bay 7 (c.1-c.3) - see Fig. \ref{fig:meshNineBay} for bay numbering. 
    The RB of the nonlinear ROMs consists of 15 VMs and 35 SMDs. The nonlinear ROMs are in good agreement with the Abaqus HFM, accurately capturing smearing and shifting of resonance peaks, out of band response, as well as the increased axial displacement level with respect to the linear solution. Notice that, for symmetry reasons, the in plane displacements in a.1, a.2 and b.1 are close to zero. Thus, the lack of superposition of the ROM and HFM curves has no practical relevance. }
    \label{fig:psdNineBay}
\end{figure}
\subsubsection{Model construction time}
We report in this subsection the ROM construction time for the nine-bay panel, drawing a direct comparison between ROM tensor identification using standard EED and EED-ECSW.
Computational times for the two methods are listed in table \ref{compTimeNineBay}, following the same classification described in section
\ref{sec:constrTimeCurvedPanel}.
Tensor identification based on EED-ECSW allowed to reduced the total ROM construction time from 13.759 h to 0.993 h, corresponding to a speed up of 13.86.
This computational advantage comes to a large extent from the acceleration of the operations performed by Abaqus and in reading of the tangent stiffness matrices, which were respectively 21.17 and 23.63 times faster when the reduced mesh was used. 
\begin{table}[h!]
\centering
\vspace{0.2cm}
\textbf{Reduction basis construction}\\
\vspace{0.1cm}
\begin{tabularx}{0.8\textwidth}{>{\raggedright\arraybackslash}X 
   >{\raggedright\arraybackslash}X 
   >{\raggedright\arraybackslash}X
   >{\raggedright\arraybackslash}X }
   \hline
    VMs & 132 s &\   & \ \\
    %\hline
    SMDs & 584 s &\   &\  \\
    \hline
    \textbf{total RB} & $\bm{716}$\ \textbf{s} & \  &\  \\
\end{tabularx}\\
\vspace{0.2cm}
\textbf{Tensor identification}\\
\vspace{0.1cm}
\begin{tabularx}{0.8\textwidth}{>{\raggedright\arraybackslash}X 
   >{\raggedright\arraybackslash}X 
   >{\raggedright\arraybackslash}X
   >{\raggedright\arraybackslash}X }
  \hline
    &EED  & EED-ECSW  & speed-up\\
    \hline
    reduced mesh & - & 0.238 h & -\\
    %\hline
    Abaqus & 8.533 h &  0.318 h & 26.83 \\
    %\hline
    reading & 4.914 h & 0.212 h   & 23.18\\
    %\hline
    identification & 78 s & 75.9 s& 1.03 \\
    %\hline
    transformation & 0.5 s & 0.4 s &1.25 \\
    other  & 259 s & 13 s & 19.92\\
    %\hline
    \hline
    \textbf{total Id.} & $\bm{13.56}$\ \textbf{h} & $\bm{0.794}$\ \textbf{h} &$\bm{17.09}$\\
    \hline  \hline
    \textbf{total RB \& Id.} & $\bm{13.759}$\ \textbf{h} & $\bm{0.993}$\ \textbf{h} &$\bm{13.86}$\\
    \hline
\end{tabularx}
\caption{Tensors construction time for the nine-bay panel: comparison between EED and EED-ECSW. All operations were run on Euler cluster of ETH Z\"urich, using 10 cores with 3Gb RAM each.}
\label{compTimeNineBay}
\end{table}
\subsubsection{Time integration and results}
The ROMs with tensors constructed by EED and EED-ECSW were integrated in Matlab using the Newmark-$\beta$ integration scheme. 
On the other hand, time integration of the HFM was performed in Abaqus using the HHT scheme.
A constant time step of $1.429 \cdot 10^{-4} \ \text{s}$ was adopted for all the integration routines, and $10 \ \text{s}$ real time were simulated.
For comparison purposes, a ROM of the linearized model was constructed with a RB composed of the first 90 VMs (to cover all the excitation bandwidth), and run for the same pressure load applied to the nonlinear models.
Displacements time histories were recorded for nine different nodes on the skin panel, each of them in the center of a different bay.
During post processing, the PSD of displacements history ($70,002$ samples each) was estimated using the Welch's method, with segments of $4,000$ data points and $50 \%$ overlap, leading to a frequency resolution of $1.75\ \text{Hz}$. In Fig. \ref{fig:psdNineBay}, we present the  PSD for displacements degrees of freedom corresponding to nodes at the center of bays 5,7 and 8 (see Fig. \ref{fig:meshNineBay}).\\
Power spectral density of displacements in the out of plane direction (along $z$ direction) is well approximated by the two ROMs, as can be seen from Fig. \ref{fig:psdNineBay} a.3, b.3, c.3.
Resonance peaks are lowered, shifted to higher frequencies and smeared out, as compared to the linear model, suggesting hardening behavior.\\
Even more remarkable differences between linear and nonlinear response are visible in the PSD of in plane displacements (see Fig. \ref{fig:psdNineBay} a.1, b.1, c.1 for displacements along $x$; see Fig. \ref{fig:psdNineBay} a.2, b.2, c.2 for displacements along $y$), which are way larger for the nonlinear models.
The excitation of axial dofs is triggered by geometric nonlinearities, responsible for the bending stretching coupling in flat structures.
Again, departure from the linear regime is well captured by the ROMs for all  monitored nodes, except for the axial dof at center of bay 5 (Fig \ref{fig:psdNineBay}.a.2).
However, it should be noticed that since the panel and the load distribution are symmetric with respect to the $x$ axis passing through the center of bay 5, axial displacements along $y$ are almost zero. 
Notice that similar arguments can be made to justify the low nodal displacements level represented in Figs. \ref{fig:psdNineBay}.a.1, \ref{fig:psdNineBay}.b.1.\\ 
From a ROM perspective, it is relevant to notice that the efficient tensor identification proposed in this work, does not affect the quality of the ROM approximation, as the PSD curves for the two ROMs are perfectly overlaid.
As such, lack of fidelity of the ROM in reproducing the HFM solution, has to be attributed to ROM model construction choices, independent of the hyper-reduction approximation introduced for tensor construction (e.g. choice of the RB).\\
Good accuracy and limited computational times make the ROM the best solution for vibration assessment.
The ROM with tensors constructed using EED-ECSW was integrated in only 20 minutes, achieving a speed up of 267 with respect to Abaqus simulation that required 89.12 h to complete.
A similar speed up was observed for the ROM with tensors constructed using EED-ECSW, as expected.
\subsubsection{Effect of training ECSW tolerance on accuracy and speed-up}
\begin{figure}
    \centering
    \includegraphics[width=1\linewidth]{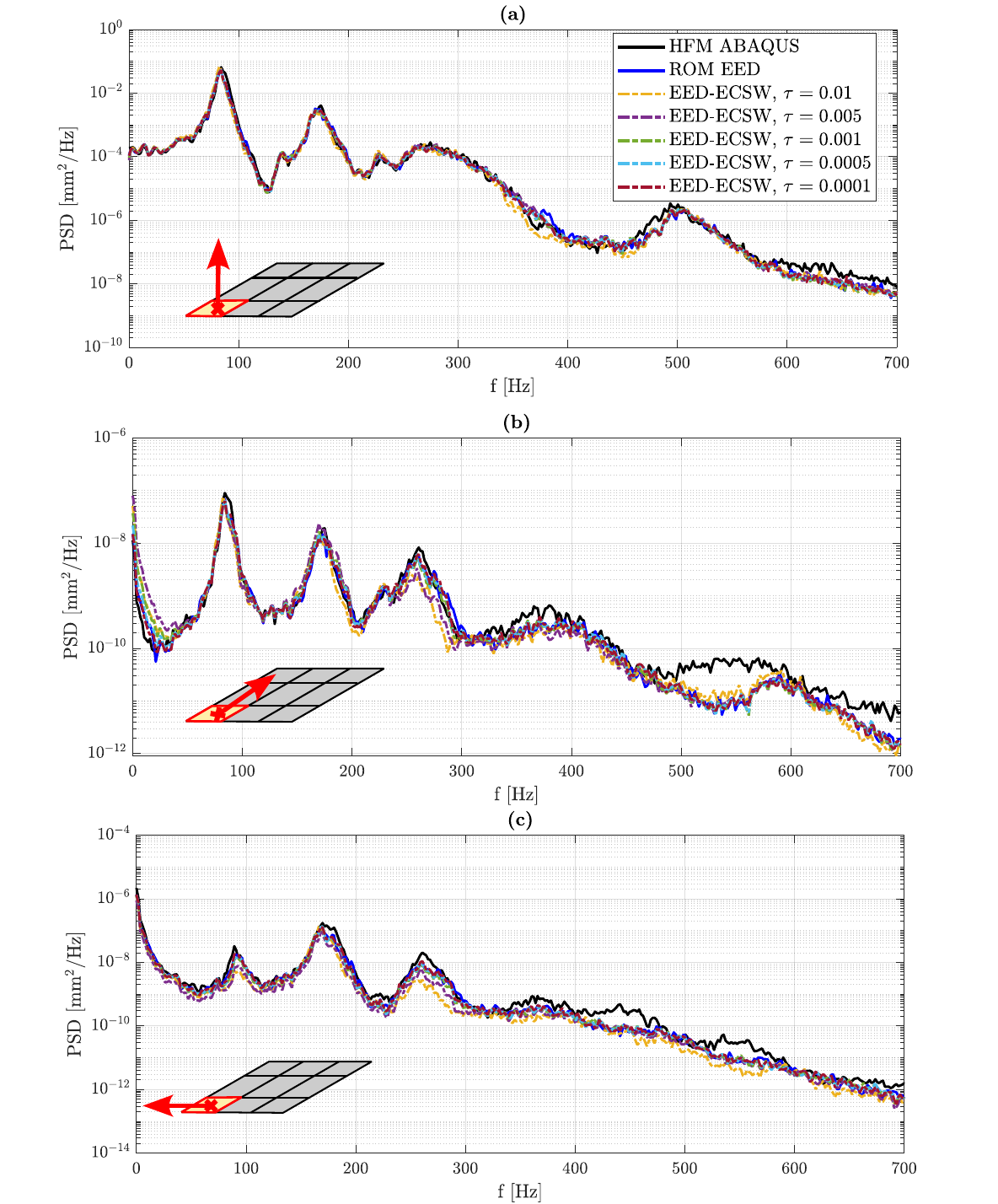}
    \caption{Effect of tolerance $\tau$ used in ECSW training on PSD accuracy of displacements for the nine-bay panel. Different reduced tensors are identified with EED-ECSW using $\tau = 0.01,0.005,0.001,0.0005,0.0001$ for the computation of the reduced mesh. The solutions obtained with tensors identified with EED and  with the HFM are shown for comparison purposes. Displacements along $z$ direction, $x$ direction and $y$ direction for a node in the middle of bay $7$ are plotted respectively in (a), (b) and (c). Marginal loss in accuracy is shown in the axial displacements (see (b) and (c)) as $\tau$ increases, while out of plane displacements (a) are less affected.}
    \label{fig:ECSWtolNineBay}
\end{figure}

The choice of the relative tolerance value $\tau$ to use in ECSW training is driven by experience of the analyst.
In this subsection we illustrate how the relative tolerance impacts on the performance of the ROM, both in terms of construction efficiency and prediction accuracy.
To this end, we identified ROM tensors with EED-ECSW using different values of relative tolerance $\tau = 0.01,0.005,0.001,0.0005,0.0001$.
As expected, we obtained different reduced meshes consisting of different number of elements, affecting the reduced mesh validation error $\epsilon_{\text{ECSW}}$ defined in \eqref{eq:ecswErr} and the computational time for the ensuing tensor identification.
The result of this experiment are reported in table \ref{tensIdTable}. 
The speed up in tensor identification obtainable using EED-ECSW instead of EED is inversely proportional to values of the relative tolerance $\tau$.
\begin{table}[h!]
\centering
\begin{tabularx}{0.8\textwidth}{>{\raggedright\arraybackslash}X 
   >{\raggedright\arraybackslash}X 
   >{\raggedright\arraybackslash}X
   >{\raggedright\arraybackslash}X
   >{\raggedright\arraybackslash}X }
   \hline
    $\tau$ & $\%\ \nelr$ & $\% \  \epsilon_{\text{ECSW}}$   & $t_c\ [s]$ & speed-up \\
    \hline
    0.01 & 0.37 &\  1.71 &\ 2000& 24.4 \\ 
    0.005 & 0.56 &\   0.94&\ 2179& 22.40\\
    0.001 & 1.19 &\   0.25&\ 2857& 17.09\\
    0.0005 & 1.58 &\   0.15&\ 3263& 14.96\\
    0.0001 & 2.89 &\   0.05&\ 4621& 10.56\\
    \hline
    \end{tabularx}\\
\caption{Effect of relative tolerance $\tau$ used in ECSW model training on the percentage of elements retained in the reduced mesh ($\%\ \nelr$), on the percentage ECSW validation error $\% \  \epsilon_{\text{ECSW}}$, on the total time for tensor identification with EED-ECSW ($t_c$) and on the speed-up in tensor identification with respect to EED. All operations were run on Euler cluster of ETH Z\"urich, using 10 cores with 3Gb RAM each.}
\label{tensIdTable}
\end{table}
Conversely, accuracy decreases as $\tau$ increases, as demonstrated by the rise in the validation error $\epsilon_{\text{ECSW}}$ for larger values of $\tau$.
To better quantify accuracy of the identified tensors we computed the PSDs of displacements using the ROMs corresponding to the different $\tau$ values, and compared them with the displacements obtained using the exact tensors, as shown in Fig.\ref{fig:ECSWtolNineBay}.
As can be seen from the figure, an increase in $\tau$ leads to a decrease in accuracy, with the axial displacements along $x$ and $y$ directions (see Fig.\ref{fig:ECSWtolNineBay}.b,Fig.\ref{fig:ECSWtolNineBay}.c) more affected than the out of plane displacements along $z$ direction (see Fig.\ref{fig:ECSWtolNineBay}.a).
%Notice however that the for the axial dof at the center of bay 5 (Fig 9.a.2) the ROM seems to be inaccurate

%% Conclusion 
\section{Conclusion} \label{sec:Conclusion}
We introduced a new procedure to accelerate the EED scheme for non-intrusive identification of nonlinear stiffness tensors in projection based ROM with RB of VMs and MDs.
The speed up is achieved by replacing the evaluations of the reduced tangent stiffness matrix with fast-to-compute approximations based on ECSW.
Reduced mesh and associated weights are obtained using a simulation-free approach where training displacements snapshots are generated from quasi-random samples of the parametrization coordinates of the SQM.
Reduced Order Model tensors are firstly identified in a physical RB using EED and subsequently transformed to their counterparts associated to an orthogonalized RB.
This step is essential to avoid imposing non-physical displacements to the FE model.
\\
The proposed methodology was tested on two different case studies: a shallow rectangular curved panel and a reinforced nine-bay panel subjected to random acoustic loading.
In both these applications, the broad-bandwidth nature of the load triggers many structural modes, resulting in large ROMs whose tensors are expensive to construct.
Dynamic response was investigated using a ROM with tensors identified using standard EED, the same ROM with approximated tensors identified using EED-ECSW, the HFM and a linear ROM.
Power spectral densities of displacements computed from simulations performed with ROM equipped with approximated tensors are in excellent agreement with the ones computed using the exact tensor ROM and the HFM.
Nonlinear effects, such as smearing and shifting of resonance peaks and out-of-band response are  well captured by the nonlinear models (ROMs and HFM), whereas they are missed by the linearized ROM. \\
The proposed approach for efficient tensor identification allowed for a remarkable reduction of the computational time for ROM construction when compared to standard EED.
Speed ups of 3.18 and 13.86 in ROM construction were recorded for the curved panel and the nine-bay panel respectively.
The provided test cases show that the proposed method is appealing for the design of industrial structural components in which different ROMs have to be constructed for different tested geometries, requiring therefore efficient ROM construction.

\section{Appendix A}

\label{sec:appendix}
We here prove that the SMPFs introduced in Eq. \eqref{eq:MPF}
are the linear combination coefficients of  the linear static solution to load vector $\mathbf{p}$, for  unit normalized modes.
Starting from the e.o.m. for a linear mechanical conservative system
\begin{equation}
     \M \qdd + \K \q= \mathbf{p}
\end{equation}
we expand the solution in a basis of VMs
\begin{equation}
    \q = \sum_{i=1}^{n} \vm_i \qrs_i
    \label{eq:provExp}
\end{equation}
and project the e.o.m. on the modes obtaining $n$ decoupled scalar equations in the modal coordinates $\eta_i$, obtaining
\begin{equation}
    \vmT_i \M \vm_i \qrdds_i + \vmT_i \mathbf{K} \vm_i \qrs_i = \vmT_i \mathbf{p},\quad \text{for} \ i \in \{1,...,n\}.
\end{equation}
These equations can be compactly written as 
\begin{equation}
    \qrdds_i + \omega_i^2 \qrs_i = \frac{1}{\mu_i}\vmT_i \mathbf{p} ,\quad \text{for} \ i \in \{1,...,n\},
\end{equation}
where $\mu_i = \vmT_i \M \vm_i$ and $\omega_i$ are respectively the modal mass and angular frequency of mode $i$.
By plugging the static solution to load vector $\mathbf{p}$  
\begin{equation}
    \qrs_i = \frac{1}{\mu_i \omega_i^2}\vmT_i \mathbf{p} = \frac{1}{\vm_i^T\mathbf{K}\vm_i}\vmT_i \mathbf{p}  
\end{equation}
in the expansion \eqref{eq:provExp} we get
\begin{equation}
    \q = \sum_{i=1}^{n} \frac{\vm_i}{\|\vm_i\|_2} \left( \frac{\|\vm_i\|_2}{\vm_i^T\mathbf{K}\vm_i}\vmT_i \mathbf{p} \right) = \sum_{i=1}^{n} \frac{\vm_i}{\|\vm_i\|_2} \text{sMPF}_i 
\end{equation}
where we unit normalized the VMs in the RB, by dividing and multiplying vectors $\vm_i$ by their 2 norm.

\bibliography{RRC2023_bib}

\begin{thebibliography}{10}
\expandafter\ifx\csname url\endcsname\relax
  \def\url#1{\texttt{#1}}\fi
\expandafter\ifx\csname urlprefix\endcsname\relax\def\urlprefix{URL }\fi
\expandafter\ifx\csname href\endcsname\relax
  \def\href#1#2{#2} \def\path#1{#1}\fi

\bibitem{Hollkamp2018}
J.~J. Hollkamp, Experiences with nonlinear modeling and acoustic fatigue, Journal of Sound and Vibration 437 (2018) 437--446.
\newblock \href {https://doi.org/10.1016/j.jsv.2018.04.029} {\path{doi:10.1016/j.jsv.2018.04.029}}.

\bibitem{Gordon2011}
R.~W. Gordon, J.~J. Hollkamp, Reduced-order models for acoustic response prediction of a curved panel, 2011.
\newblock \href {https://doi.org/10.2514/6.2011-2081} {\path{doi:10.2514/6.2011-2081}}.

\bibitem{Spottswood2008}
S.~Spottswood, J.~Hollkamp, T.~Eason, On the use of reduced-order models for a shallow curved beam under combined loading, American Institute of Aeronautics and Astronautics (AIAA), 2008.
\newblock \href {https://doi.org/10.2514/6.2008-2235} {\path{doi:10.2514/6.2008-2235}}.

\bibitem{Przekop2007}
A.~Przekop, S.~A. Rizzi, Dynamic snap-through of thin-walled structures by a reduced-order method, AIAA Journal 45 (2007) 2510--2519.
\newblock \href {https://doi.org/10.2514/1.26351} {\path{doi:10.2514/1.26351}}.

\bibitem{przekop2006nonlinear}
A.~Przekop, S.~A. Rizzi, D.~S. Groen, Nonlinear acoustic response of an aircraft fuselage sidewall structure by a reduced-order analysis, in: Ninth International Conference on Recent Advances in Structural Dynamics, no. Paper-135, 2006.

\bibitem{mignolet2013review}
M.~P. Mignolet, A.~Przekop, S.~A. Rizzi, S.~M. Spottswood, A review of indirect/non-intrusive reduced order modeling of nonlinear geometric structures, Journal of Sound and Vibration 332~(10) (2013) 2437--2460.

\bibitem{krysl2001dimensional}
P.~Krysl, S.~Lall, J.~E. Marsden, Dimensional model reduction in non-linear finite element dynamics of solids and structures, International Journal for numerical methods in engineering 51~(4) (2001) 479--504.

\bibitem{kerschen2007physical}
G.~Kerschen, F.~Poncelet, J.-C. Golinval, Physical interpretation of independent component analysis in structural dynamics, Mechanical Systems and Signal Processing 21~(4) (2007) 1561--1575.

\bibitem{lieu2005pod}
T.~Lieu, C.~Farhat, M.~Lesoinne, Pod-based aeroelastic analysis of a complete f-16 configuration: Rom adaptation and demonstration, in: 46th AIAA/ASME/ASCE/AHS/ASC structures, structural dynamics and materials conference, 2005, p. 2295.

\bibitem{geradin2015mechanical}
M.~G{\'e}radin, D.~J. Rixen, Mechanical vibrations: theory and application to structural dynamics, John Wiley \& Sons, 2015.

\bibitem{Idelsohn1985}
S.~R. Idelsohn, A.~Cardona, A reduction method for nonlinear structural dynamic analysis, Computer Methods in Applied Mechanics and Engineering 49 (1985) 253--279.
\newblock \href {https://doi.org/10.1016/0045-7825(85)90125-2} {\path{doi:10.1016/0045-7825(85)90125-2}}.

\bibitem{idelsohn1985load}
S.~R. Idelsohn, A.~Cardona, A load-dependent basis for reduced nonlinear structural dynamics, Computers \& Structures 20~(1-3) (1985) 203--210.

\bibitem{Marconi2020}
J.~Marconi, P.~Tiso, F.~Braghin, A nonlinear reduced order model with parametrized shape defects, Computer Methods in Applied Mechanics and Engineering 360 (3 2020).
\newblock \href {https://doi.org/10.1016/j.cma.2019.112785} {\path{doi:10.1016/j.cma.2019.112785}}.

\bibitem{Marconi2021}
J.~Marconi, P.~Tiso, D.~E. Quadrelli, F.~Braghin, A higher-order parametric nonlinear reduced-order model for imperfect structures using neumann expansion, Nonlinear Dynamics 104 (2021) 3039--3063.
\newblock \href {https://doi.org/10.1007/s11071-021-06496-y} {\path{doi:10.1007/s11071-021-06496-y}}.

\bibitem{Morteza}
M.~K. Mahdiabadi, P.~Tiso, A.~Brandt, D.~J. Rixen, A non-intrusive model-order reduction of geometrically nonlinear structural dynamics using modal derivatives, Mechanical Systems and Signal Processing 147 (1 2021).
\newblock \href {https://doi.org/10.1016/j.ymssp.2020.107126} {\path{doi:10.1016/j.ymssp.2020.107126}}.

\bibitem{Sombroek2018}
C.~S. Sombroek, P.~Tiso, L.~Renson, G.~Kerschen, Numerical computation of nonlinear normal modes in a modal derivative subspace, Computers and Structures 195 (2018) 34--46.
\newblock \href {https://doi.org/10.1016/j.compstruc.2017.08.016} {\path{doi:10.1016/j.compstruc.2017.08.016}}.

\bibitem{Wu2016}
L.~Wu, P.~Tiso, Nonlinear model order reduction for flexible multibody dynamics: a modal derivatives approach, Multibody System Dynamics 36 (2016) 405--425.
\newblock \href {https://doi.org/10.1007/s11044-015-9476-5} {\path{doi:10.1007/s11044-015-9476-5}}.

\bibitem{Rutzmoser2017}
J.~B. Rutzmoser, D.~J. Rixen, A lean and efficient snapshot generation technique for the hyper-reduction of nonlinear structural dynamics, Computer Methods in Applied Mechanics and Engineering 325 (2017) 330--349.
\newblock \href {https://doi.org/10.1016/j.cma.2017.06.009} {\path{doi:10.1016/j.cma.2017.06.009}}.

\bibitem{tiso2011reduction}
P.~Tiso, E.~Jansen, M.~Abdalla, Reduction method for finite element nonlinear dynamic analysis of shells, AIAA journal 49~(10) (2011) 2295--2304.

\bibitem{Farhat2014}
C.~Farhat, P.~Avery, T.~Chapman, J.~Cortial, Dimensional reduction of nonlinear finite element dynamic models with finite rotations and energy-based mesh sampling and weighting for computational efficiency, International Journal for Numerical Methods in Engineering 98 (2014) 625--662.
\newblock \href {https://doi.org/10.1002/nme.4668} {\path{doi:10.1002/nme.4668}}.

\bibitem{Jain2018}
S.~Jain, P.~Tiso, Simulation-free hyper-reduction for geometrically nonlinear structural dynamics: A quadratic manifold lifting approach, Journal of Computational and Nonlinear Dynamics 13 (7 2018).
\newblock \href {https://doi.org/10.1115/1.4040021} {\path{doi:10.1115/1.4040021}}.

\bibitem{An2008}
S.~S. An, T.~Kim, D.~L. James, Optimizing cubature for efficient integration of subspace deformations, ACM Transactions on Graphics 27 (12 2008).
\newblock \href {https://doi.org/10.1145/1409060.1409118} {\path{doi:10.1145/1409060.1409118}}.

\bibitem{rutzmoserThesis}
J.~Rutzmoser, Model order reduction for nonlinear structural dynamics, Ph.D. thesis, Technische Universit\"at M\"unchen (2018).

\bibitem{Mignolet2008}
M.~P. Mignolet, C.~Soize, Stochastic reduced order models for uncertain geometrically nonlinear dynamical systems, Computer Methods in Applied Mechanics and Engineering 197 (2008) 3951--3963.
\newblock \href {https://doi.org/10.1016/j.cma.2008.03.032} {\path{doi:10.1016/j.cma.2008.03.032}}.

\bibitem{nashThesis}
M.~Nash, Nonlinear structural dynamics by finite elemetn modal synthesis, Ph.D. thesis, Imperial College London (1977).

\bibitem{shi1996finite}
Y.~Shi, C.~Mei, A finite element time domain modal formulation for large amplitude free vibrations of beams and plates, Journal of sound and vibration 193~(2) (1996) 453--464.

\bibitem{almroth1978automatic}
B.~O. Almroth, P.~Stern, F.~A. Brogan, Automatic choice of global shape functions in structural analysis, AIAA Journal 16~(5) (1978) 525--528.

\bibitem{Muravyov2003}
A.~A. Muravyov, S.~A. Rizzi, Determination of nonlinear stiffness with application to random vibration of geometrically nonlinear structures, Computers and Structures 81 (2003) 1513--1523.
\newblock \href {https://doi.org/10.1016/S0045-7949(03)00145-7} {\path{doi:10.1016/S0045-7949(03)00145-7}}.

\bibitem{McEwan2001}
M.~I. McEwan, J.~R. Wright, J.~E. Cooper, A.~Y. Leung, A finite element/modal technique for nonlinear plate and stiffened panel response prediction, American Institute of Aeronautics and Astronautics Inc., 2001.
\newblock \href {https://doi.org/10.2514/6.2001-1595} {\path{doi:10.2514/6.2001-1595}}.

\bibitem{Hollkamp2008}
J.~J. Hollkamp, R.~W. Gordon, Reduced-order models for nonlinear response prediction: Implicit condensation and expansion, Journal of Sound and Vibration 318 (2008) 1139--1153.
\newblock \href {https://doi.org/10.1016/j.jsv.2008.04.035} {\path{doi:10.1016/j.jsv.2008.04.035}}.

\bibitem{Perez2014}
R.~Perez, X.~Q. Wang, M.~P. Mignolet, Nonintrusive structural dynamic reduced order modeling for large deformations: Enhancements for complex structures, Journal of Computational and Nonlinear Dynamics 9 (7 2014).
\newblock \href {https://doi.org/10.1115/1.4026155} {\path{doi:10.1115/1.4026155}}.

\bibitem{Wang2018}
X.~Q. Wang, G.~P. Phlipot, R.~A. Perez, M.~P. Mignolet, Locally enhanced reduced order modeling for the nonlinear geometric response of structures with defects, International Journal of Non-Linear Mechanics 101 (2018) 1--7.
\newblock \href {https://doi.org/10.1016/j.ijnonlinmec.2018.01.007} {\path{doi:10.1016/j.ijnonlinmec.2018.01.007}}.

\bibitem{Perez2011}
R.~Perez, X.~Q. Wang, M.~P. Mignolet, Reduced order modeling for the nonlinear geometric response of cracked panels, 2011.
\newblock \href {https://doi.org/10.2514/6.2011-2018} {\path{doi:10.2514/6.2011-2018}}.

\bibitem{farhat2020}
C.~Farhat, S.~Grimberg, A.~Manzoni, A.~Quarteroni, et~al., Computational bottlenecks for proms: precomputation and hyperreduction, Model order reduction 2 (2020) 181--244.

\bibitem{chaturantabut2010nonlinear}
S.~Chaturantabut, D.~C. Sorensen, Nonlinear model reduction via discrete empirical interpolation, SIAM Journal on Scientific Computing 32~(5) (2010) 2737--2764.

\bibitem{Farhat2015}
C.~Farhat, T.~Chapman, P.~Avery, Structure-preserving, stability, and accuracy properties of the energy-conserving sampling and weighting method for the hyper reduction of nonlinear finite element dynamic models, International Journal for Numerical Methods in Engineering 102 (2015) 1077--1110.
\newblock \href {https://doi.org/10.1002/nme.4820} {\path{doi:10.1002/nme.4820}}.

\bibitem{Trainotti2024}
F.~Trainotti, J.~Marinko, J.~Maierhofer, D.~J. Rixen, Ecsw hyperreduction of hyper-viscoelastic components via co-simulation with abaqus, Finite Elements in Analysis and Design 241 (11 2024).
\newblock \href {https://doi.org/10.1016/j.finel.2024.104222} {\path{doi:10.1016/j.finel.2024.104222}}.

\bibitem{Kim2023}
Y.~Kim, S.~H. Kang, H.~Cho, H.~Kim, S.~J. Shin, Parametric reduced-order modeling enhancement for a geometrically imperfect component via hyper-reduction, Computer Methods in Applied Mechanics and Engineering 403 (2023) 115701.
\newblock \href {https://doi.org/10.1016/J.CMA.2022.115701} {\path{doi:10.1016/J.CMA.2022.115701}}.

\bibitem{crisfield1991}
M.~Crisfield, Nonlinear Finite Element Analysis of Solids and Structures, Wiley, 1991.

\bibitem{tiso20214}
P.~Tiso, M.~Karamooz~Mahdiabadi, J.~Marconi, Modal methods for reduced order modeling, in: Model Order Reduction: Volume 1: System-and Data-Driven Methods and Algorithms, De Gruyter, 2021, pp. 97--138.

\bibitem{Andersson2023}
L.~Andersson, P.~Persson, K.~Persson, Efficient nonlinear reduced order modeling for dynamic analysis of flat structures, Mechanical Systems and Signal Processing 191 (5 2023).
\newblock \href {https://doi.org/10.1016/j.ymssp.2023.110143} {\path{doi:10.1016/j.ymssp.2023.110143}}.

\bibitem{Jain2017}
S.~Jain, P.~Tiso, J.~B. Rutzmoser, D.~J. Rixen, A quadratic manifold for model order reduction of nonlinear structural dynamics, Computers and Structures 188 (2017) 80--94.
\newblock \href {https://doi.org/10.1016/j.compstruc.2017.04.005} {\path{doi:10.1016/j.compstruc.2017.04.005}}.

\bibitem{Weeger2016}
O.~Weeger, U.~Wever, B.~Simeon, On the use of modal derivatives for nonlinear model order reduction, International Journal for Numerical Methods in Engineering 108 (2016) 1579--1602.
\newblock \href {https://doi.org/10.1002/nme.5267} {\path{doi:10.1002/nme.5267}}.

\bibitem{tiso2011optimal}
P.~Tiso, Optimal second order reduction basis selection for nonlinear transient analysis, in: Modal Analysis Topics, Volume 3: Proceedings of the 29th IMAC, A Conference on Structural Dynamics, 2011, Springer, 2011, pp. 27--39.

\bibitem{Rutzmoser2017QM}
J.~B. Rutzmoser, D.~J. Rixen, P.~Tiso, S.~Jain, Generalization of quadratic manifolds for reduced order modeling of nonlinear structural dynamics, Computers and Structures 192 (2017) 196--209.
\newblock \href {https://doi.org/10.1016/j.compstruc.2017.06.003} {\path{doi:10.1016/j.compstruc.2017.06.003}}.

\bibitem{Vizzaccaro2021}
A.~Vizzaccaro, L.~Salles, C.~Touz\'{e}, Comparison of nonlinear mappings for reduced-order modelling of vibrating structures: normal form theory and quadratic manifold method with modal derivatives, Nonlinear Dynamics 103 (2021) 3335--3370.
\newblock \href {https://doi.org/10.1007/s11071-020-05813-1} {\path{doi:10.1007/s11071-020-05813-1}}.

\bibitem{Veraszto2020}
Z.~Veraszto, S.~Ponsioen, G.~Haller, Explicit third-order model reduction formulas for general nonlinear mechanical systems, Journal of Sound and Vibration 468 (3 2020).
\newblock \href {https://doi.org/10.1016/j.jsv.2019.115039} {\path{doi:10.1016/j.jsv.2019.115039}}.

\bibitem{Lin2023}
J.~Lin, X.~Q. Wang, B.~Wainwright, M.~P. Mignolet, Improved identification of stiffness coefficients of non intrusive nonlinear geometric reduced order models of structures, International Journal of Non-Linear Mechanics 152 (6 2023).
\newblock \href {https://doi.org/10.1016/j.ijnonlinmec.2023.104380} {\path{doi:10.1016/j.ijnonlinmec.2023.104380}}.

\bibitem{helton2003latin}
J.~C. Helton, F.~J. Davis, Latin hypercube sampling and the propagation of uncertainty in analyses of complex systems, Reliability Engineering \& System Safety 81~(1) (2003) 23--69.

\bibitem{Schoneman2017}
J.~D. Schoneman, M.~S. Allen, R.~J. Kuether, Relationships between nonlinear normal modes and response to random inputs, Mechanical Systems and Signal Processing 84 (2017) 184--199.
\newblock \href {https://doi.org/10.1016/j.ymssp.2016.07.010} {\path{doi:10.1016/j.ymssp.2016.07.010}}.

\bibitem{hayes1996statistical}
M.~H. Hayes, Statistical Digital Signal Processing and Modeling, John Wiley \& Sons, New York, 1996.

\bibitem{buehrle2000finite}
R.~D. Buehrle, G.~A. Fleming, R.~S. Pappa, F.~W. Grosveld, Finite element model development for aircraft fuselage structures, in: XVIII International Modal Analysis Conference, 2000.

\end{thebibliography}

\end{document}